\pdfoutput=1

\documentclass[preprint,11pt]{elsarticle}
\usepackage[letterpaper,top=1in, bottom=1in, left=1in, right=1in]{geometry}
\usepackage{graphicx}
\usepackage{amssymb}
\usepackage{lineno}
\usepackage[hidelinks]{hyperref}
\usepackage{amsmath, amsthm}
\usepackage{float}

\theoremstyle{definition}

\usepackage{algorithm,algcompatible}
\usepackage{algpseudocode}
\usepackage{algorithmicx}
\usepackage{booktabs}
\usepackage{subfig}
\usepackage{multirow}
\usepackage{lscape}
\usepackage{makecell}
\usepackage[english]{babel}
\usepackage[markup=underlined]{changes}
\definechangesauthor[color=red]{JLM}
\definechangesauthor[color=blue]{JH}

\usepackage{todonotes}
\setcommentmarkup{\todo[color={authorcolor!20},size=\scriptsize]{#3: #1}}

\newcolumntype{M}[1]{>{\centering\arraybackslash}m{#1}}
\DeclareMathAlphabet{\mathpzc}{OT1}{pzc}{m}{it}
\usetikzlibrary{positioning}

\journal{}

\begin{document}

\begin{frontmatter}


\title{Charging-as-a-Service: On-demand battery delivery for light-duty electric vehicles for mobility service}

\author[ua_address]{Shuocheng Guo}
\ead{sguo18@ua.edu}

\author[ua_address]{Xinwu Qian\corref{cor}}
\cortext[cor]{Corresponding author.}
\ead{xinwu.qian@ua.edu}

\author[ua_address]{Jun Liu}
\ead{jliu@eng.ua.edu}

\address[ua_address]{Department of Civil, Construction and Environmental Engineering, The University of Alabama, Tuscaloosa, AL 35487, United States}

\begin{abstract}

This study presents an innovative solution for powering electric vehicles, named Charging-as-a-Service (CaaS), that concerns the potential large-scale adoption of light-duty electric vehicles (LDEV) in the Mobility-as-a-Service (MaaS) industry. Analogous to the MaaS, the core idea of the CaaS is to dispatch service vehicles (SVs) that carry modular battery units (MBUs) to provide LDEVs for mobility service with on-demand battery delivery. The CaaS system is expected to tackle major bottlenecks of a large-scale LDEV adoption in the MaaS industry due to the lack of charging infrastructure and excess waiting and charging time. A hybrid agent-based simulation model (HABM) is developed to model the dynamics of the CaaS system with SV agents, and a trip-based stationary charging probability distribution is introduced to simulate the generation of charging demand for LDEVs. Two dispatching algorithms are further developed to support the optimal operation of the CaaS. The model is validated by assuming electrifying all 13,000 yellow taxis in New York City (NYC) that follow the same daily trip patterns. Multiple scenarios are analyzed under various SV fleet sizes and dispatching strategies. The results suggest that an optimal CaaS system with 250 SVs may serve the LDEV fleet in NYC with an average waiting time of 5 minutes, save the travel distance to visit charging station at over 50 miles per minute, and gain considerable profits of up to $\$50$ per minute. This study offers significant insights into the feasibility, service efficiency, and financial sustainability for deploying city-wide CaaS systems to power the electric MaaS industry.

\end{abstract}


\end{frontmatter}

\section{Introduction}

The transportation sector has witnessed two major revolutions in the past decade that have been reshaping the mobility landscape of urban travelers: one being the rapid growing adoption of electric vehicles (EVs) in replace of the traditional gasoline vehicles advances, while the other attributing to the unprecedented rise of the Mobility-as-a-Services (MaaS). In particular, over 2.1 million EVs were sold in 2019 globally, marking a 40\% annual growth rate since 2010~\cite{outlook2020entering}. Meanwhile, the MaaS market has seen tremendous growth since 2014, with the number of ride-hailing trips increased from 60,000 to over 750,000 in New York City (NYC) alone. Given the vast number of daily trips from the MaaS market, the use of EVs signifies opportunities to alleviate the worse emission and energy consumption situations in large cities~\cite{qian2020impact}. Nevertheless, the two streams of revolution barely converge. Despite the ambitious plan from the major transportation network companies (TNCs) for electrifying their fleet~\cite{CleanTrans}, less than 0.2\% of the fleet for Uber and Lyft are EVs.

There are at least three critical barriers that prevent the widespread commercial adoption of light-duty electric vehicles (LDEVs) in the MaaS industry. The first barrier attributes to the lack of sufficient charging facilities in urban areas. For instance, there were 450 public charging stations in NYC providing 931 Level 2 chargers and 82 DC Fast chargers as of 2019~\cite{doecharge}, which can hardly satisfy the enormous charging needs if we fully electrify the fleet of major MaaS operators with over 100,000 vehicles. Tesla has recently announced the Supercharger fair use policy~\cite{tesla_ban}, banning commercial EVs from using its Supercharger stations. The second barrier arises from insufficient land space and power grid supply to allow for the massive construction of charging facilities in MaaS's core service areas. This indicates that the charging facilities have to be located at distant locations, which leads to the third barrier due to excessive travel time to-and-from a charging station, long charging time, and prolonged waiting time at the station during peak hours. During peak hours in Shenzhen, China, electric taxi drivers have to wait at least 30 minutes for an available charger\cite{dong2017rec}. These barriers pose emerging needs for a dedicated charging solution that tailors to MaaS in order to promote the electrification of the MaaS industry.

Existing studies of charging facility planning primarily focused on locating fixed charging stations (FCS), battery swapping stations (BSS), and more recently, the deployment of mobile charging systems (MCS). Early work proposed the optimal design of FCS locations under the network equilibrium of the coupled road and power networks~\cite{he2013optimal}, recharging time and flow-dependent energy~\cite{he2014network}, and the multi-class user equilibrium~\cite{he2015deploying}. Moreover, the concept of charging lanes envisions wireless charging and motivates studies on the deployment of the charging lanes~\cite{riemann2015optimal,chen2016optimalb} and its co-planning problem with existing FCSs~\cite{chen2017deployment}. More recently, researchers investigated the optimal design of both location and capacity for FCSs under the consideration of long-distance travel~\cite{wang2019designing}, electric bus charging stations~\cite{he2019fast,lin2019multistage}, the waiting and charging spots for electric taxis~\cite{yang2017data}, and road users with travel and charging equilibrium~\cite{chen2020optimal}. 
 
For the deployment of BSSs, the objectives of the early work included developing a cost-effective subscription based swapping service~\cite{mak2013infrastructure}, maximizing the net present value or profit~\cite{zheng2013electric,sarker2014optimal}, minimizing the long-term operating costs~\cite{sun2014optimal}, and reducing queuing time under a city-wide electric taxi fleet~\cite{wang2017toward}. While there exist plentiful studies on FCS and BSS, the research on mobile charging facilities is still in the early stage. The early study considered a mobile charging platform consisting of portable plug-in chargers and mobile swapping stations~\cite{huang2014design}. With the incorporation of fast charging technology, more recent studies explored the operation of commercial mobile charging facilities including mobile charging stations~\cite{chen2016optimala} and mobile swapping vans~\cite{shao2017mobile}.

Although the aforementioned studies on charging facility preparation can promote the adoption of privately owned EVs, they may hardly benefit LDEVs for the MaaS industry considering the distinct operation dynamics. We briefly summarize the features of the charging facilities, including their barriers to promote the MaaS industry in Table~\ref{tab:summary_of_features}. As discussed previously, these key barriers consist of the high implementation cost, limited space available for charging facilities, the interruption of mobility services for charging batteries, and the difficulties in coordinating the charging needs among large-scale LDEVs. We note that the mobile modular units, such as prototypes of SparkCharge~\cite{SparkCharge}, promises a more flexible and portable scheme. The mobile unit consists of one modular charger unit and up to five modular battery units (MBUs). The MBU can stack on top of one another and work together, with each MBU offers at least 12 miles of range under the charging speed of at least one mile per minute~\cite{SparkCharge}, which is sufficient to satisfy the service needs of LDEVs in urban areas and makes ``charge-as-you-go" possible. And the above barriers can be largely resolved if there exists a system that can optimally coordinate the delivery of the MBUs with the charging needs from the LDEVs. Aside from the daily use scenarios, we remark that the CaaS using the MBUs will also be more resilient to extreme scenarios  (e.g., charging demand surge in peak hours, power outage, severe weather), outperforming the mobile charging station due to its portability and flexibility. With a sufficiently large SV fleet and additional MBUs as the back-up, the CaaS can ensure the functionality of the LDEV fleet under those extreme cases to maintain a high service level and constant power supply. These features motivate us to explore the innovative charging solution for powering LDEV fleet and investigate the feasibility and performances of the city-wide adoption of such charging systems.

\begin{table}[H]
\small
\centering
\caption{Summary of charging facilities and infrastructure}
\begin{tabular}{p{1.5cm}p{1.5cm}p{1.5cm}p{3cm}p{4cm}l}
\toprule
Facility type                      & Charging speed (mile/min)         & Location                      & Configuration (per unit)                                                                          & Barriers for the MaaS                      & Source              \\ \hline
Fixed charger                    & 0.03 - 4 & Fixed &   Installation costs:\$4,000-\$51,000(DC fast)                                                                              & High cost, space limitation, long charging time, service interruption for charging, difficult to cooperate large-scale charging events & \cite{FCSchargingspeed}~\cite{smith2015costs} \\ \hline
Battery swapping station                           & Swap in 2.5 min   & Fixed              &       Construction cost: \$800,000; maintenance costs: \$30,000/year                                                                                   & High cost, space limitation, service interruption for charging, restriction on the type of battery used & \cite{wang2017toward}~\cite{lidicker2011business}    \\ \hline
Mobile charging station & Up to 0.5         & Flexible                      & Weight: 1860 lbs; battery capacity: 80 kWh                                               & Not portable, capacity and weight limitation for delivery &\cite{MobiEVcharger}  \\ \hline
Mobile modular units         & 1              & Flexible                      & Weight: 19.8 lbs (MCU), 48.4 lbs (MBU);  battery capacity: 3.5 kWh & Small capacity, requires the support of optimal delivery algorithms   &\cite{SparkCharge} \\ 
\bottomrule
\end{tabular}
\label{tab:summary_of_features}
\end{table}

In this study, we make the initial attempt to investigate an emerging charging solution for electrifying the MaaS industry, named Charging-as-a-Service (CaaS). Unlike the current notion of CaaS, which mainly offers subscription service for fixed chargers~\cite{EVPassport}~\cite{EVgo}, the CaaS discussed in this study is a more typical on-demand charging service via MBU delivery which is motivated by the recent advances in the modular fast-charging technology~\cite{SparkCharge,GMCharge}. We consider \textit{a centralized CaaS system that dispatches service vehicles (SVs) to the LDEV's requested charging location to replace the depleted MBUs with fully-charged ones}. A hybrid agent-based model(HABM) is developed to simulate the joint operation of LDEVs and the SV fleet in a CaaS system, aiming at providing LDEVs an efficient charging service while evaluating the operation performances, benefits, and financial sustainability of developed CaaS systems. The main contributions of our work can be summarized as follows:

\begin{enumerate}
    \item A CaaS operation framework is designed to provide on-demand battery delivery services to satisfy the changing needs of city-wide LDEV fleet for mobility service. 
    \item A hybrid agent-based model is developed which consists of stationary charging demand generation for large-scale LDEV fleet and agent-level service dynamics for the SV fleet. 
    \item Dispatching algorithms are developed to support the optimal operation of SVs that minimize the average waiting time for the charging requests. 
    \item Real-world numerical experiments are conducted to examine the applicability of the CaaS system, which examines the service dynamics and quantifies the savings of out-of-service time and charging miles and its financial sustainability when compared to charging at the FCSs.
\end{enumerate}

The rest of the paper is organized as follows. The next section introduces the business model of CaaS, the HABM framework, and the performance metrics. In Section 3, we present the numerical experiments for electrifying the NYC taxi market as a case study and discuss the results and insights. Finally, section 4 concludes our work with key findings, recommendations of CaaS system configurations, and future directions.

\section{Model}

\subsection{Preliminaries}

This study models the dynamics of the CaaS system at the zonal level with discrete time steps. In the CaaS system, there are two types of interacting dynamics: (1) the movement of LDEVs, both serving passenger trips as well as relocations, that results in charging requests (demand) and (2) the movement of SVs that provide the MBU delivery services to the site of service requests (supply). The LDEVs will make reservations for the MBU replacement before running out of electricity and meet SVs at the scheduled location and time. The out-of-electricity LDEVs will resume their trip service upon the delivery of the MBUs, with MBUs function as power banks that provide sufficient power to sustain travels in urban areas (1 mile per minute). The SVs will keep serving charging requests from LDEVs until consuming all fully-charged MBUs, and the SVs will then return to the depot to unload empty MBUs and reload with fully charged ones. The interactions between the LDEV and SV fleet are ruled by the dispatching policies that assign the available SVs to the requested LDEVs to minimize the LDEV fleet's total waiting time.

\begin{figure}[H]
    \centering
    \includegraphics[width=0.75\textwidth]{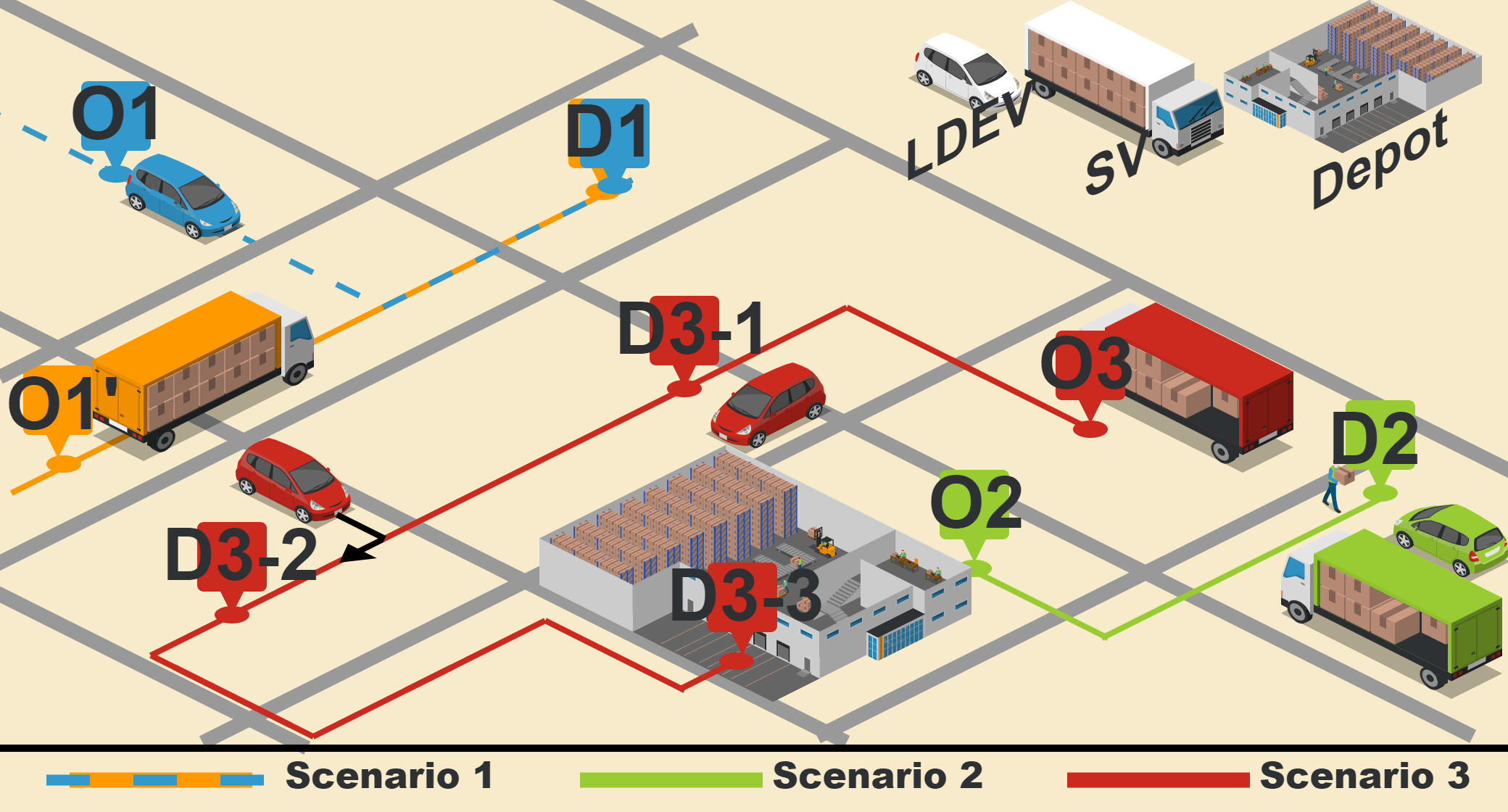}
    \caption{An illustration of the CaaS system with three basic scenarios consisting of LDEV, SV, and the depot. This illustration chart is modified from an open-source template in the Icograms~\cite{Icograms}.}
    \label{fig:illustration}
\end{figure}

We start with three basic scenarios that illustrate the operation of the CaaS system, as shown in Figure~\ref{fig:illustration}. Scenario 1 shows how the CaaS system processes reservation requests for charging service. The blue LDEV from the origin (O1) sends out a request with its destination (D1) and estimated arrival time to D1. The system then matches this request with the orange SV at O1$^{'}$, and the LDEV and SV will meet at the D1. Scenario 2 represents an LDEV being served, in which both the green LDEV and the assigned green SV have reached the destination (D2), and the SV driver is replacing the depleted MBUs with fully-charged ones. Scenario 3 describes that the red SV has two scheduled red LDEVs with requested destinations at D3-1 and D3-2. When the SV is en route, one LDEV is waiting for SV at D3-1 and the other is on its way to D3-2. When the SV is out of MBUs after serving the request at D3-2, it will move back to the depot at D3-3 to replenishment fully-charged MBUs. The following sections will introduce the details on the generation of charging demand, the operation of the CaaS in the hybrid agent-based model, and three different dispatching policies.

\subsection{Generation of stationary charging demand}\label{sec:generation of charging demand}

One bottleneck for planning charging infrastructure for the LDEVs is the ability to characterize the distribution of charging needs. Unlike private EVs where charging locations can be easily estimated, the location of charging needs of LDEVs depending on the stochastic trip and relocation activities. To address this issue, we introduce a stochastic model which \textit{considers moving activities in the MaaS system as a discrete-time Markov chain}, following the preliminary analyses from the authors' prior work~\cite{qian2019stationary}. With the zonal representation, the Markov chain regards each zone as a state and the movement between zone $i$ and zone $j$ as the transition probability $P_{ij}$. In this regard, the MaaS system dynamics are modeled as an event-based system with each event being either serving a trip or performing a relocation. The sequences of moving events will drive the LDEVs between zones and eventually incur charging requests. 

Specifically, we consider a centralized MaaS service system where the LDEVs follow trip and relocation sequences designated by the platform. We divide the operation hours into discrete time intervals, where the passenger demand pattern is assumed to be stable within each time interval. While a specific trip sequence may differ for individual LDEVs, their collective movement patterns at time $t$ can be measured statistically. That is, the probability $\pi_i^t(k)$ that an arbitrary LDEV will remain at zone $i$ following $k$th event can be described as:
\begin{equation}
    \pi_i^t(k)=\sum_{j=1}^N \pi_j^t(k-1)P_{ji}^t
\end{equation}
With a sufficient large number of events (trips), we will have the stationary probability that an arbitrary LDEV will be in a specific zone following any event as:
\begin{equation}
    \bar{\pi}^t=\bar{\pi}^t P^t
\end{equation}
Moreover, for each time interval $t$, the battery of an LDEV at the end of each movement event may fall into one of the two categories: (1) sufficient battery and continue to the next movement or (2) recharge the battery by the CaaS and then continue to the next movement. Denote $E_{ij}$ as the expected electricity consumption between zone $i$ and zone $j$, we can therefore represent the expected battery consumption for any LDEV that arrives in zone $i$ as: 
\begin{equation}
\bar{e}_i^t =\frac{\sum_{j=1}^{N} \bar{\pi}_{j}^tP_{ji}^tE_{ji}^t}{\sum_{j=1}^{N}
\bar{\pi}_{j}^{t}P_{ji}^t}
\label{eq:stationary_electricity}
\end{equation}
We note that both $P_{ij}^t$ and $E_{ij}^t$ can be obtained by mining the historical daily operation data from the large-scale LDEV fleet. Let $C_{LDEV}$ be the battery capacity and $\delta$ be the charging threshold (e.g., an LDEV will not charge until the battery level is below $\delta C_{LDEV}$, the probability that an arbitrary LDEV that arrives in zone $i$ and has its battery below the threshold (therefore needs recharge) follows:

\begin{equation}
    P^t(X_{i}^{charge}=1)=\frac{\bar{\pi}_i^t\bar{e}_i^t}{C_{LDEV}(1-\delta)}
    \label{eq:charge_chance}
\end{equation}
where a detailed derivation can be found in~\cite{qian2019stationary}.
Let $K$ be the average number of trips that an LDEV may complete per unit time, $N_{LDEV}$ be the LDEV fleet size and $N_w$ be the number of LDEVs that are awaiting recharging of the battery, we can eventually express the arrival rate of charging demand per unit time in a given zone $i$ as:
\begin{equation}
    \lambda_i^t=K\bar{\pi}_i^t(N_{LDEV}-N_{w})P^t(X_{i}^{charge}=1)
    \label{eq:charge_demand}
\end{equation}
where we have $N_{LDEV}-N_{w}$ be the set of active LDEVs that will generate charging requests, and $K\bar{\pi}^t_i$ represents the expected number of trips that may arrive in zone $i$ per unit time per LDEV. We can then use the Poisson arrival process to generate the arrival charging demand following equation~\ref{eq:charge_demand}. This helps avoid the expensive simulation of the actual moving trajectories for individual LDEVs by directly sampling the charging demand at the zonal level, which is much more scalable and contributes to significant speedups of the simulation when the fleet size is large.

\subsection{The hybrid agent-based model}

\begin{figure}[h!]
\centering
\includegraphics[width=\textwidth]{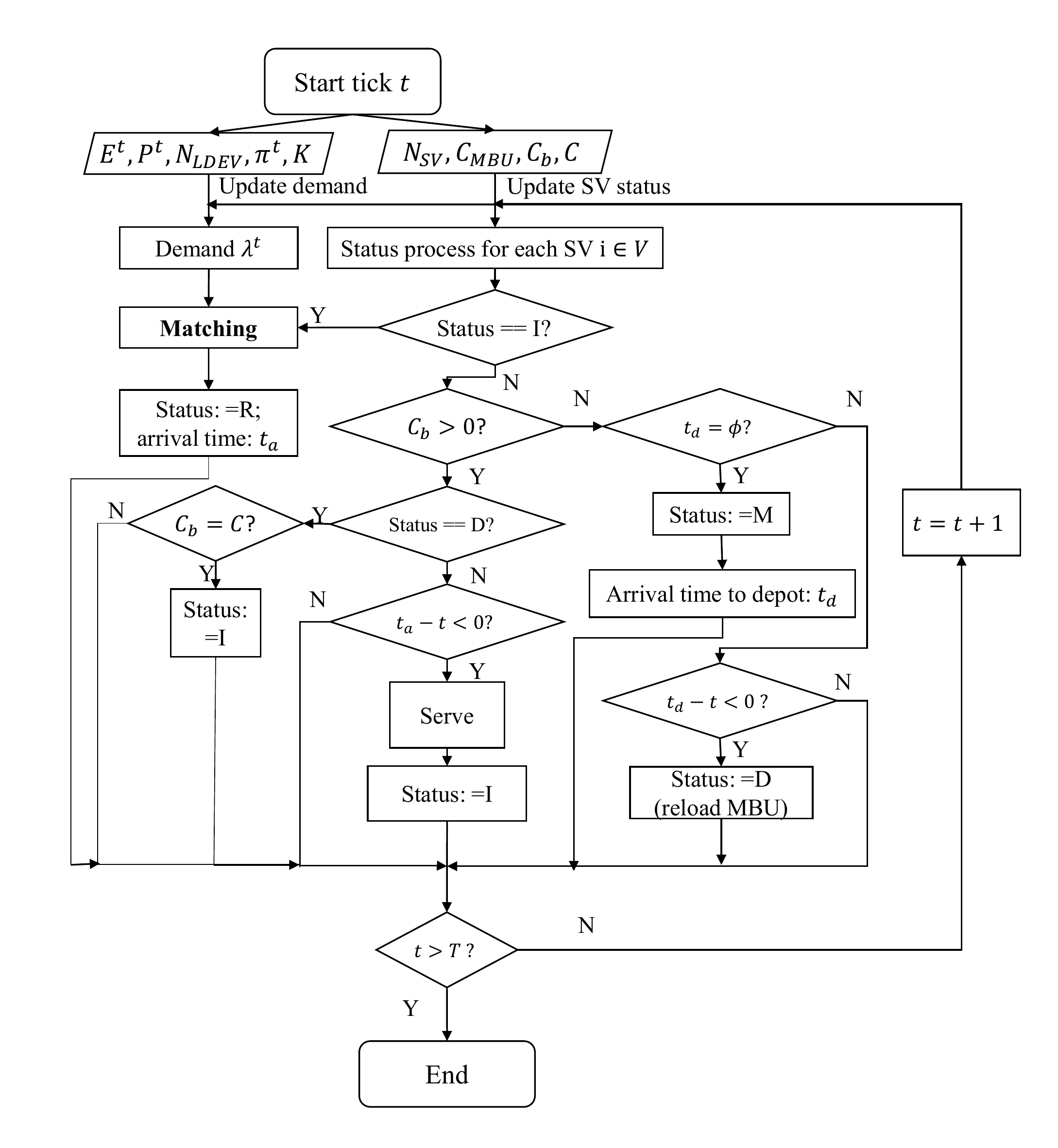}
\caption{Integrated HABM framework for modeling the CaaS system. }
\label{fig:flowchart}
\end{figure}

Given the arrival of charging demand, the dynamics for the fleet of SVs are simulated following the agent-based modeling framework, as shown in Figure~\ref{fig:flowchart}. The simulation framework takes the SV fleet size ($N_{SV}$) and battery capacity of the MBU ($C_{MBU}$) as input and follows three different policies to match available SVs to charging requests. And $C_{b}$ denotes the remaining number of fully-charged MBUs on the SV, which is limited by the capacity of the SV ($C$). The outputs are the performance metrics of the CaaS system in serving the generated charging requests that are parameterized by $E$,$P$,$N$,$\bar{\pi}$, and $K$.

The framework describes the criterion on how the SV agents are connected to the charging requests, as discussed in the preliminary section, and how the SV agents are navigated between the charging depot and charging requests in the HABM. Specifically, we define that each SV agent may be in one of the following five states:
\begin{itemize}
\item \textbf{Idle (I)}. An SV is in I state if it has sufficient fully-charged MBUs to serve future charging requests and is currently not assigned with a delivery task.
\item \textbf{Relocating (R)}. Once an SV is matched with a charging request, it will remain in R state until it arrives at the requested service location. 
\item \textbf{Ready for service (S)}. An SV is in S state if it arrives at the requested service location, and the SV will then replace the depleted MBU with a fully-charged one. 
\item \textbf{Moving back to depot (M)}. When an SV runs out of the fully-charged MBU, it will be assigned in D state and start moving back to the depot location. 
\item \textbf{At depot (D)}. An SV is in D state if it has returned to the depot. The SV will then unload the depleted MBUs and reload with fully-charged MBUs. 
\end{itemize}

To develop the HABM framework, we make the following assumptions:
\begin{enumerate}
    \item We assume that each SV will have sufficient range to deliver all the MBUs without any interruption, and the SV will restore its full range while at the depot.
    \item Since a zone-based HABM framework is adopted, the relocating distance between zone $i$ and zone$j$ is sampled from the distribution of observed travel duration (e.g., from historical data) between the two zones, which is assumed to follow a log-normal distribution in our study.
    \item The SVs are operated by a central dispatching platform. This means that the SV drivers need to follow assigned service order and will not cruise for orders on their own.
\end{enumerate}

In each simulation tick $t$, the states of SVs will be based on if there are charging requests in the service area and if there are available fully-charged MBUs. For regular MBU replacement services, an SV will undergo the I-R-S-I process as the charging requests arrive: from I to R when matched with a charging request, from R to S while relocating to the service location, and from S back to I when the service is completed. The I-R-S-I process continues until the SV runs out of fully-charged MBUs, where the SV will then be redirected to follow the S-M-D-I process to send the depleted MBUs back to the depot and reload with fully-charged MBUs to maximum capacity. The SV will change from S to M when the service is completed, and all fully-charged MBUs are consumed, from M to D when the SV is navigated back to the depot for reloading new MBUs, and from D to I when the SV completed the reloading process and is again available for the services. 

\subsection{Dispatching policy}
One key component of the HABM framework, as shown in Figure~\ref{fig:flowchart}, is the matching module, which sets up the rule on how available SVs will be dispatched to satisfy the charging requests. In this study, we consider three dispatching policies that match available SVs with charging requests. The first policy, as a benchmark scenario, is the greedy first-come-first-serve policy (FCFS). The other two policies, namely the optimal matching with local information (OML) and optimal matching with dynamic detour (OMDD), are developed under the consideration of minimizing total matching cost (defined as the waiting time of charging requests). They differ in the set of SVs that are available for the dispatching policy at each time step. OML can only see SVs that are currently in idle state, while OMDD can control the entire SV fleet despite of the SVs' state. In this regard, OML is more suitable for the CaaS system with independent SV drivers while OMDD is designed for fully-cooperative fleet. 

\subsubsection{First-come-first-serve}
Under the FCFS policy, LDEVs' charging requests are served in the order in which they are received with the nearest available SV. The charging requests are sorted in the ascending order in the request list $\mathcal{L}$ based on their request time. At each simulation tick, requests are matched sequentially with available SVs, and travel cost $\mathcal{T}_{ij}$ is incurred for a request in zone $i$ that is matched with an SV in zone $j$. We note that $\mathcal{T}_{ij}$ is a random variable to reflect the impacts of stochastic locations of the charging request and the LDEV in their corresponding zones. For each simulation tick, the FCFS terminates when either all charging requests are matched or all available SVs are occupied. All matched charging requests are then removed from $\mathcal{L}$. Under the FCFS policy, the matching can be solved with time complexity of $\mathcal{O}(|\mathcal{L}|log(|\mathcal{L}|)+|\mathcal{L}|N)$, where $N$ denotes the number of available SVs.

We note that the FCFS policy will likely lead to more vacant trips and result in low SV fleet efficiency with its sequential matching scheme. For instance, when charging requests emerge between the central area and distant area in turns, the SVs have to travel a long distance back and forth to satisfy all requests. In this case, all LDEVs have to wait long for the SVs, and the low level of service potentially results in a higher unsatisfied rate. Next, we will explore two optimal dispatching policies that improve both LDEV users' level of service and SV fleet efficiency.

\subsubsection{Optimal matching with local information} 

OML policy aims to minimize the total waiting cost of charging requests, including relocation cost and current waiting cost. With local information, the OML policy will only dispatch SVs that are currently in I state. In each tick $t$, candidate SVs are assigned following an optimal dispatching protocol $f (\mathcal{V}^{t},\mathcal{D}^{t})$, where $\mathcal{V}^{t}$ denotes the set of available SVs and $\mathcal{D}^{t}$ for the set of LDEVs that need MBU replacements. When LDEV $k$ in zone $i$ schedules a MBU replacement at destination $j$ in time $t$, the platform will set its future service time as $T_i=\tau_{ij}^t+t$ for the charging request, where $\tau_{ij}^{t}$ is the estimated relocation time from zone $i$ to zone $j$ in time $t$. 

To derive an optimal assignment between SVs and LDEVs, we formulate the OML strategy as a minimum weight bipartite matching problem with $\mathcal{V}^{t},\mathcal{D}^{t}$ as the two sets of nodes for each time step $t$. Considering edges connecting each pair of nodes in the two sets, the weight on each edge is measured as the difference between $c_{ij}^{t}+t$ and $T_j$, for all $i\in \mathcal{V}^{k}$ and $j\in \mathcal{D}^{k}$. With $c_{ij}^t$ being the relocation time from SV $i$ to request $j$, the time difference may lead to early arrival of SV if $c_{ij}^t+t$ is smaller than $T_j$, or late arrival otherwise. In our study, we penalize both early and late arrivals but with a penalty factor $0\leq \eta\leq 1$ to favor early arrival over late arrival. In this regard, the OML dispatching policy at time $t$ can be mathematically formulated as:

\begin{equation}
\begin{aligned}
    &\min_{x} \, \sum_{i \in \mathcal{V}^{t}} \sum_{j \in \mathcal{D}^{t}} x_{ij} \left(\max(0,c_{ij}^t + t-T_j) +\eta \max(0, T_j-c_{ij}^t-t)\right)\\
    \text{subject to}\quad  & \sum_{i \in \mathcal{V}^{t}} x_{ij}=1, \forall\, i\in\mathcal{V}^{t}\\
    &\sum_{j \in \mathcal{D}^{t}} x_{ij}=1, \forall\, i\in\mathcal{D}^{t}\\
    &x_{ij}\in\{0,1\}\\
\end{aligned}
\label{eq:match}
\end{equation}

The constraints in equation~\ref{eq:match} state that each available SV is assigned to only one LDEV, and vice versa. In most cases, the cardinality of $\mathcal{V}^t$ and $\mathcal{D}^t$ will not be equal so that a perfect matching can not be established. To address this issue, we augment the set of smaller cardinality with $||\mathcal{V}|-|\mathcal{D}||$ dummy elements, and mark the weight on edges that connect to dummy modes with an arbitrarily large value $M$. Consequently, one can solve equation~\ref{eq:match} with augmented sets to find the optimal dispatching using the Hungarian method~\cite{kuhn1955hungarian}, which can solve the problem in $\mathcal{O}(N^3)$ time ($N$ being the cardinality of the augmented set).  

As compared to the FCFS policy, the OML performs matching at the system level, which delivers optimal local dispatching of SVs. However, the assignment using local information may not ensure long-term optimal dispatching, as new requests being generated after the matching decision has been made. This motivates the third policy as an optimal matching with a dynamic detour to tackle this issue.

\subsubsection{Optimal matching with dynamic detour}

OMDD seeks to mitigate the OML policy's drawback due to its lack of flexibility for dynamically updating the assignment with the arrival of new charging requests. This is especially the case when an SV is scheduled to arrive much earlier than the service time, and the lead time is sufficient to allow for additional services. 

The OMDD inherits the overall framework of OML, with the major differences being that (1) an SV can have more than one service request scheduled, e.g., the schedule of SV $j$ being $S_j=\{m_j^1,m_j^2,...,m_j^n\}$, and (2) the matching cost is derived by inspecting if a valid gap is available for additional service in between successive scheduled requests. At time step $t$, we define that an additional service $m_j^{x}=(x,T_{x})$ can be inserted in between SV $j$'s consecutive requests $m_j^{k},m_j^{k+1}$, if the following conditions are satisfied:

\begin{equation}
    T_{k+1}-T_{k}\geq \tau_{kx}^t + \tau_{xk+1}^t
\end{equation}
so that there are sufficient time for the SV $j$ to travel from request $k$ to $x$ and then from $x$ to $k+1$. 

In addition, the new request can be inserted to the beginning of $S_j$ if 

\begin{equation}
    T_1-t\geq \tau_{jx}^t+\tau_{x1}^t
\end{equation}
with $\tau_{jx}^t$ being the travel time from SV $j$'s current location to the location of new request $x$. Finally, the new request can always be appended to the end of $S_j$ (n+1 th position). 

Let $c_{jx}^k$ be the valid cost for inserting $x$ to the $k$th request in $S_j$, we have  

\begin{equation}
    c^k_{jx}=
    \begin{cases}
    \eta(T_1-t-\tau_{jx}^t-\tau_{x1}^t), \text{if }k=1\\
    \eta(T_{k+1}-T_{k}-\tau_{kx}^t-\tau_{xk+1}^t), \text{if } 1<k\leq n\\
    \max (0,t+\tau_{nx}^t-T_x)+\eta (\max(0,T_x-\tau_{nx}^t-t)), \text{if }k=n+1
    \end{cases}
\end{equation}

And the matching cost for SV $j$ and service request $x$ is determined by the minimum of all valid costs: 

\begin{equation}
    c^*_{jx}=\min_k c_{jx}^k
\end{equation}

With the definition of $c^*_{jx}$,  the dispatching strategy of SVs under the OMDD policy can be again obtained by solving the minimum weight bipartite matching problem for each time step $t$. Different from the OML policy, $\mathcal{V}^t$ now includes all SVs in the CaaS system, and the problem can be formulated as:

\begin{equation}
\begin{aligned}
    &\min_{x_{ij}} \, \sum_{i \in \mathcal{V}^{t}} \sum_{j \in \mathcal{D}^{t}} x_{ij} c^*_{ji} \\
    \text{subject to}\quad  
    & \sum_i x_{ij} = 1, \forall\, i\in\mathcal{V}^t\\
    &\sum_j x_{ij} = 1, \forall\, i\in\mathcal{D}^t\\
    &x_{ij} \in \{0,1\}
\end{aligned}
\label{eq:match_pre_dispatch}
\end{equation}
Similar to the OML policy, the optimal solution under the OMDD policy can also be solved by the Hungarian method in $\mathcal{O}(N^3)$ time.

\subsection{Performance metrics}
The performance of the CaaS system will be evaluated via simulations with the HABM and the three dispatching policies, with the following major performance metrics. 

The first metric calculates the savings in out-of-service duration (min) based on the gap between the average waiting time of CaaS and the time it requires to visit and charge at an FCS:

\begin{equation}
\label{eq:savingoos}
    R_{o}  = \bar{T}_{FCS} -  \bar{T}_{CaaS}*\frac{C_{LDEV}}{C_{MBU}}
\end{equation}

where $\bar{T}_{CaaS}$ is the average waiting time in CaaS and $\bar{T}_{FCS}$ captures the average time required to visit and charge at an FCS. $C_{LDEV}$ is battery capacity for the LDEV and $C_{MBU}$ is battery capacity of MBUs. In Equation~\ref{eq:savingoos},\textit{ we use $\frac{C_{LDEV}}{C_{MBU}}$ to convert the corresponding waiting time for CaaS into the full-charge-equivalent (FCE) waiting time}, so that we account for the more number of MBU replacements to achieve the same battery level as the full charging at an FCS. 

Similarly, we can quantify the saving of charging distance by computing the gap in trip miles per unit time between the LDEVs' expected travel distance to visit the nearest FCS and SVs' operation distances:

\begin{equation}
\label{eq:savingdist}
    R_{d}  = \frac{1}{\Delta T} \left( \bar{d}_{FCS}*N_{s} -  d_{SV}*\frac{C_{LDEV}}{C_{MBU}} \right)
\end{equation}
where $\Delta T$ is the duration of the simulation process, $\bar{d}_{FCS}$ is the average travel distance (mile) to FCSs, and $d_{SV}$ represents the total miles traveled by the SV fleet. 

In addition, it is also important to understand the financial sustainability of the proposed CaaS services. This requires measuring the cost and revenue incurred during daily operations. We consider that the total CaaS cost includes the capital cost (e.g., the procurement of SVs, MBUs), labor cost for SV drivers, and the operation cost. Specifically, the labor cost is quantified by served number of charging requests and operation cost is associated with the total miles traveled. We express the total CaaS cost as:

\begin{equation}
\label{eq:totalcost}
Q_{CaaS}= N_{SV} \left( Q_{SV} + C*Q_{MBU} + d_{SV}*q \right) + Q_{L}V_{s}
\end{equation}

where $Q_{SV},\,Q_{MBU}$ are the per-unit purchase cost of SV and MBU, respectively. $C$ is the SV capacity for MBUs, $q$ is the coefficient for operation cost and $Q_L$ represents the labor cost per served request.The average operation cost is set as 57.3 cents per mile following the statistics for the P70D step van used by the UPS~\cite{lammert2012eighteen}.

To evaluate the revenue generated by the CaaS, we first calculate the willingness-to-adopt of CaaS as compared to alternate charging mode based on the average waiting time and the per-service price $p$. For simplicity, we use the binary logit function to model the charging choices~\cite{discretechoice} and measure the probability $P_{CaaS}^i$ that LDEV $i$ may adopt the CaaS as follows:

\begin{align}
\label{eq:probability}
    P_{CaaS}^i & = \frac{1}{1+ \exp{\left( U_{FCS}^i-U_{CaaS}^i \right)}} \\
    U_{CaaS}^i  & =-\left( p + \kappa W_{s}^i \right) * \frac{C_{LDEV}}{C_{MBU}}-Q_{c} \\ 
    U_{FCS}^i &  = -\kappa \bar{T}_{FCS} - Q_{c}
\end{align}

In equation~\ref{eq:probability}, $U_{CaaS}^i$ and $U_{FCS}^i$ refer to LDEV $i$'s perceived utilities for using the CaaS and the FCS, respectively. $\kappa$ is the time value for LDEV users, which is considered as $\$1/min$. $Q_{c}$ is the electricity cost for fully charging the battery on a typical LDEV. Here we consider a 50-kWh battery capacity for the LDEV and the electricity price of \$0.135/kWh, so that $Q_c$ takes the value of $6.75$. Similar to the calculation of $R_o$ and $R_d$, the FCE values are also used here to convert the cost per CaaS service into an equivalent cost for a full charge of the LDEV battery, and $W_{s}^{i}$ here represents the waiting time for LDEV $i$. The revenue and cost eventually enables us to evaluate the profit level of the CaaS under different scenarios. 

The parameter setting in this study is summarized in Table~\ref{tab:pricing assumption}. In particular, $C_{MBU}$ is assumed as $14$ kWh according to the SparkCharge configuration~\cite{SparkCharge}. And we assume the SV capacity for a 4-modular MBU count $C$ as $50$, by considering the SV's payload capacity~\cite{national2010technologies}, and we ensure that the total weight of 200 ($50\times 4$) MBUs is under the payload capacity. In this case, each SV will carry 50 sets of four-module MBU, where each set can provide the LDEVs with 14 kWh of electricity or at least 48 miles of drive range.

\begin{table}[!htbp]
\caption{Parameter setting used in this study}
\centering
\label{tab:pricing assumption}
\begin{tabular}{lll}
\toprule
Parameter            & Value      & Description                          \\ 
\hline
$Q_{SV}$         & \$40,000        & SV unit price                           \\ 
$Q_{MBU}$         & \$10,000        & 4-modular MBU unit price                           \\ 
$q$            & \$0.573 per mile             & SV operation cost                           \\
$C$         & 50 count       & SV's capacity for four-module MBU                    \\ 
$C_{LDEV}$         & 50 kWh        & LDEV battery capacity                     \\ 
$C_{MBU}$         & 14 kWh        & MBU battery capacity                     \\ 
$\kappa$           & \$1 per min     & time value                           \\
$\bar{T}_{FCS}$ & 90 min          & average waiting time and driving time to FCS \\ 
$\bar{d}_{FCS}$ & 15 mile          & average relocation distance to FCS \\ 
\bottomrule
\end{tabular}
\end{table}

\section{Results}

\subsection{Data}

We choose New York City (NYC) as the study area for our numerical experiments. We leverage NYC taxi trip data in 2013~\cite{taxidata2013} to obtain occupied and vacant trip sequences, where these trips are preprocessed and aggregated into 263 taxi zones that cover Manhattan, Brooklyn, Queens, Bronx, Staten Island, and Newark Liberty International Airport (EWR). The 2013 taxi trip data is used instead of more recent data as it represents the most recent taxi dataset in NYC that has the encrypted medallion ID available, based on which we can track complete trip sequences of individual taxis. 
This data allows us to calibrate $\bar{\pi}^t, K, P^t$ for generating the charging demand in the HBAM. 

\begin{figure}[H]
    \centering
    \subfloat[Valid zones, 260 of 263~\label{fig:target_zone}]{\includegraphics[width=0.3\textwidth]{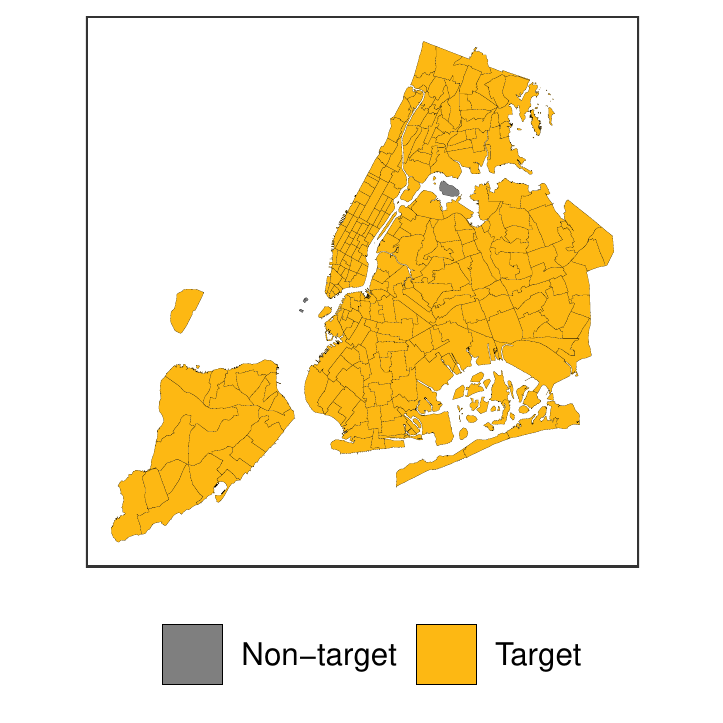}}\hfill
    \subfloat[Demand spatial distribution~\label{fig:served_EV_num}]{\includegraphics[width=0.3\textwidth]{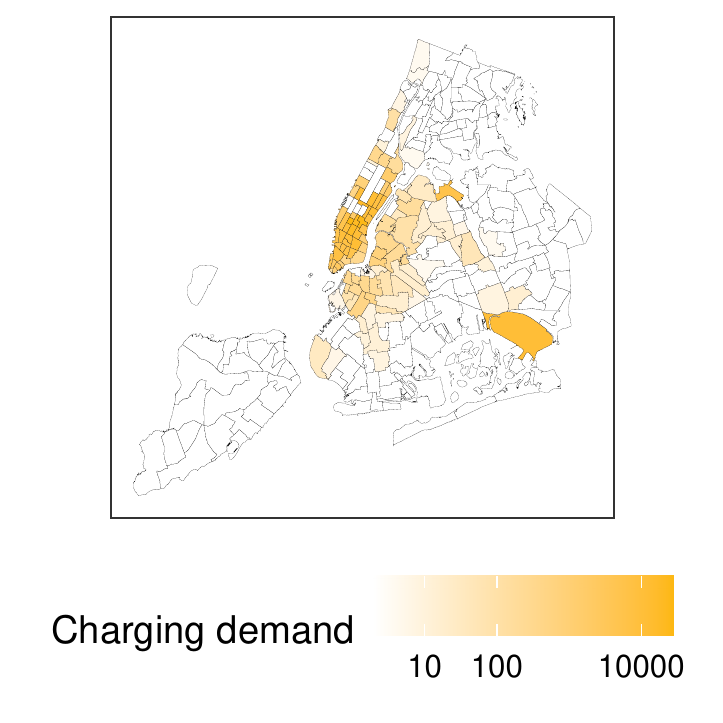}}\hfill
    \subfloat[Demand temporal distribution~\label{fig:served_demand_temporal}]{\includegraphics[width=0.3\textwidth]{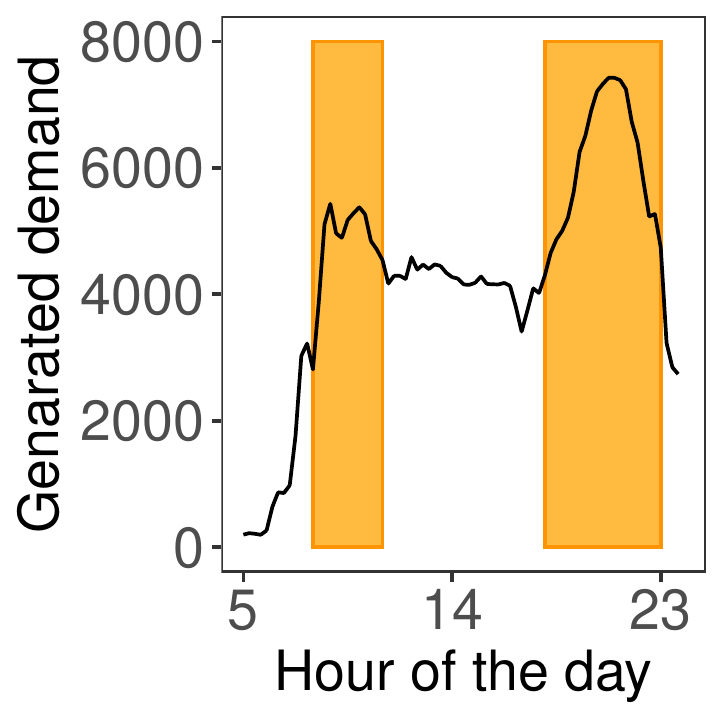}}   
    \caption{Zones with valid trips and charging demand distribution}
\end{figure}

Moreover, we prepare the distributions of the travel time between each pair of zones based on the NYC for-hire vehicles (FHV) trip record data in 2019~\cite{NYCFHV2019}. The origin-destination travel time matrix is used as the input for our HABM to determine the relocation duration and calibrate the battery consumption. We assume the average level of travel efficiency being 5 km/kWh in the NYC area~\cite{seddig2017integrating}. Among those 263 spatial zones, 260 have valid numbers of trips as shown in Figure~\ref{fig:target_zone} and these zones constitute the final study area. And we are able to measure the generated charging requests following the trip record data using equation~\ref{eq:charge_demand}. We summarize the spatial and temporal distribution of charging requests in Figure~\ref{fig:served_EV_num} and Figure~\ref{fig:served_demand_temporal}, which visualize the high-demand zones and the morning and evening peaks.

\subsection{Experiment settings}

\begin{figure}[H]
    \centering
    \subfloat[Spatial distribution of system performance for the depot candidates~\label{fig:depot_map}]{\includegraphics[width=\linewidth]{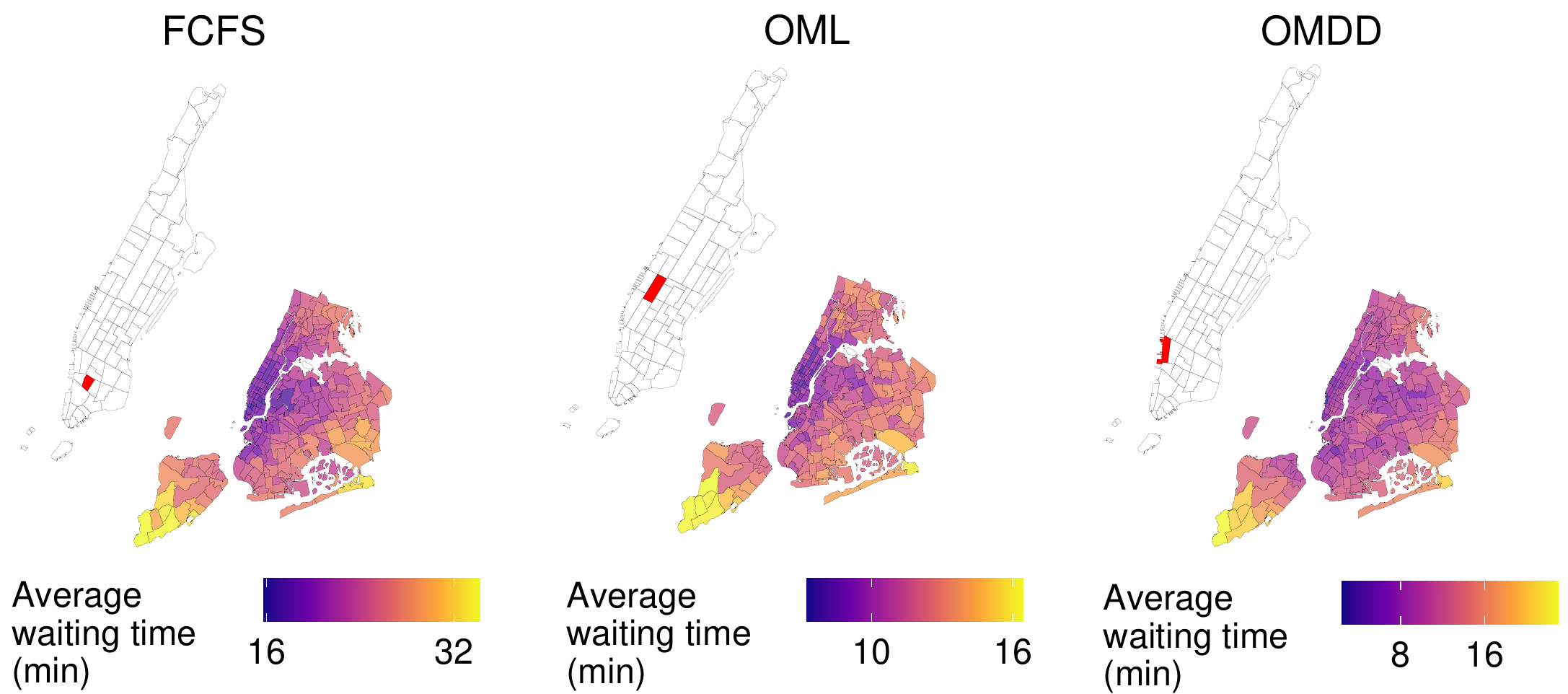}}\\
    \subfloat[Average waiting time histogram by borough (the red dashed lines represent the median)\label{fig:depot_hist}]{\includegraphics[width=.8\linewidth]{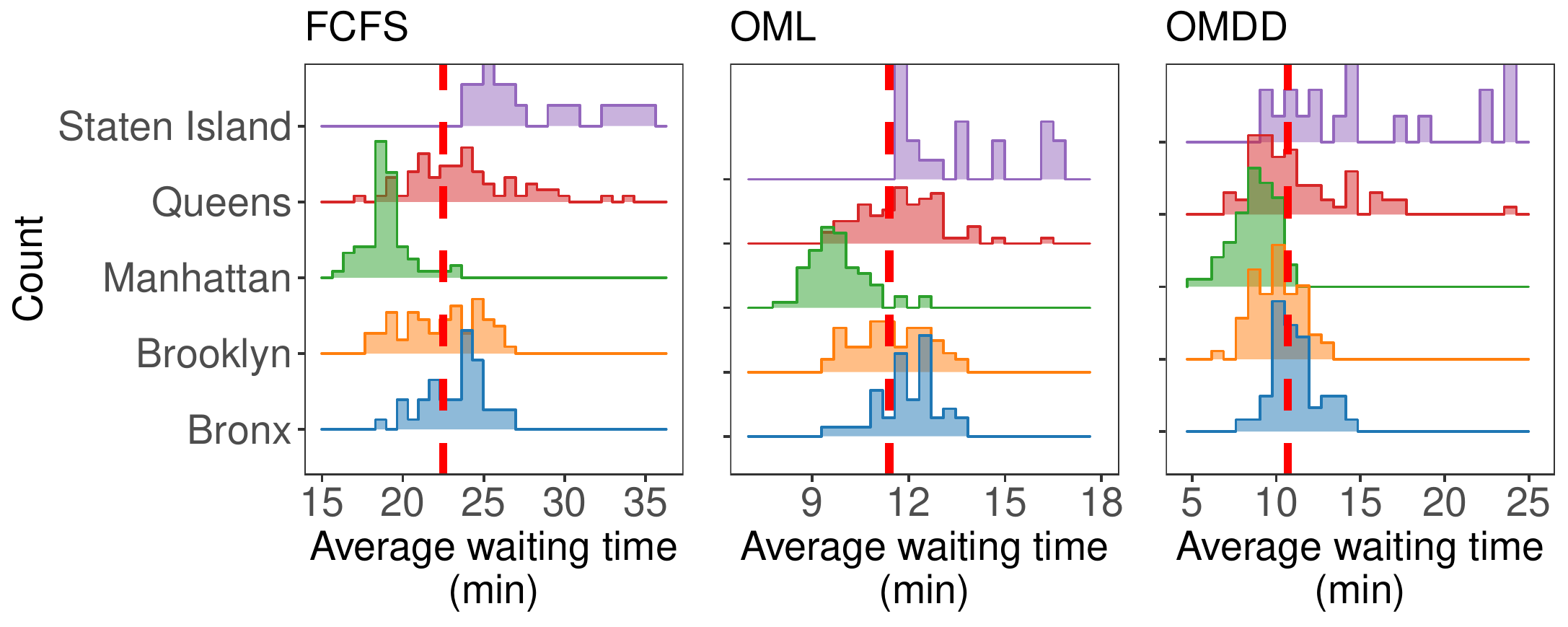}}
    \caption{The optimal depot selection process ($N_{SV}=250$)}
    \label{fig:depot_pick}
\end{figure}

We focus on the experiments for electrifying all 13,000 yellow taxis in NYC and validate the performances of CaaS for serving the charging demand of the resulting NYC LDEV fleet. We consider a single depot and set the depot to the taxi zone that achieves the minimum average waiting time. To find the optimal depot location, we assume that the SV fleet size is 250, which is large enough to satisfy most charging requests and conduct a brute-and-force searching based on the results from the HABM. And the outputs from the HABM can also be used in the simulation-optimization framework for the future extension on multi-depot locations. We pick the optimal depot locations for the three dispatching policies as shown in Figure~\ref{fig:depot_map}, which presents the distribution of the average waiting time if the depot is sited in a particular taxi zone. Unsurprisingly, the depot with the minimum average waiting time is found inside Manhattan, which is the place with most of the charging requests. In addition, Figure~\ref{fig:depot_hist} shows the histogram of average waiting time among five boroughs in NYC. The left-end bin of the histogram denotes the minimum average waiting time that can be achieved in each borough and the red dashed line marks the median value across all taxi zones.  We observe almost all zones in Manhattan have an average waiting time below the median for all three policies. However, for the OML and especially the OMDD policies, the best performing locations in Queens and Brooklyn are found to be comparable with those in Manhattan. The finding is closely related to practical concerns on siting the charging depot, as the average price per square foot commercial space is $\$684$ in Manhattan~\cite{NYClandprice2019}, but the prices are only $\$169$, $\$259$, $\$74$ in Queens, Brooklyn, and the Bronx, respectively. It indicates that the other three boroughs are much more cost-friendly and should be considered as competitive candidate locations for hosting the charging depot in light of the comparable system-level performances reported in Figure~\ref{fig:depot_hist}. Finally, the results also highlight the effectiveness of the OMDD policy since the resulting system performances are found to be less sensitive to the spatial locations of the charging facility.

After identifying the optimal depots, we run the HABM to evaluate the efficiency and profitability of the CaaS system while considering the varying SV fleet size $N$ from 125 to 500, incremented by 25, under the three dispatching strategies, the FCFS, the OML, and the OMDD. The evaluation of each scenario is performed using 20 random seeds. And we run each simulation for 1,740 iterations (ticks), with each tick representing one minute of real-world time. To get stable and representative metrics for the system performances, we run the first 600 ticks as a warm-up period and only collect results from the remaining 1,140 ticks ($ 19 \times 60 $), representing a 19-hour simulation from 5 AM to 12 AM.

\subsection{Operation performances}

This section discusses the CaaS system dynamics in detail. We first demonstrate the impacts of different SV fleet size and dispatching strategies on the CaaS system dynamics and identify one cost-effective fleet size for further analyses. And we present the system performances, temporal and spatial dynamics, and service equity in the CaaS system. Finally, we analyze the savings in travel miles and out-of-service time and the economic sustainability of the CaaS by comparing with charging at the FCSs as the alternative. 

\subsubsection{SV fleet size}

\begin{figure}[H]
\centering
    \subfloat[Average waiting time on a log scale, SD as shaded areas~\label{fig:caas_dynamics_wait}]{\includegraphics[width=0.33\linewidth]{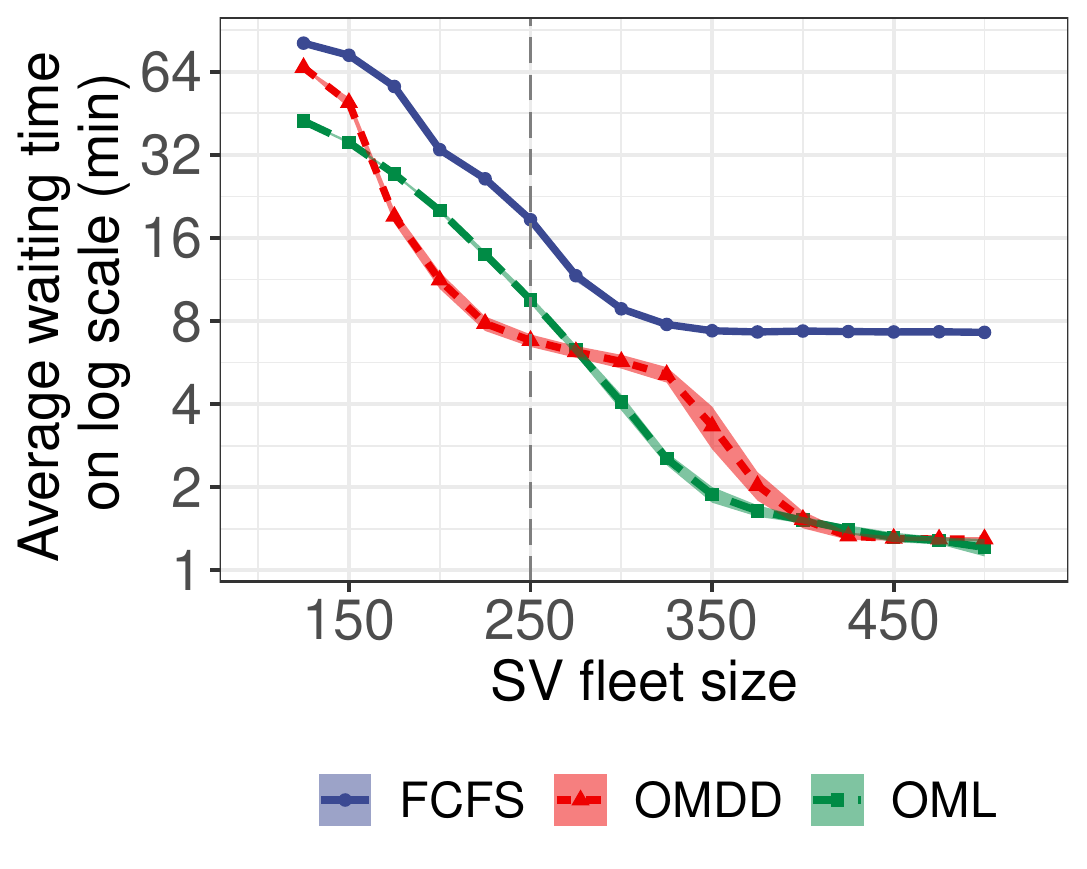}} 
    \subfloat[Number of fulfilled charging requests~\label{fig:caas_dynamics_satisfy_count}]{\includegraphics[width=0.33\linewidth]{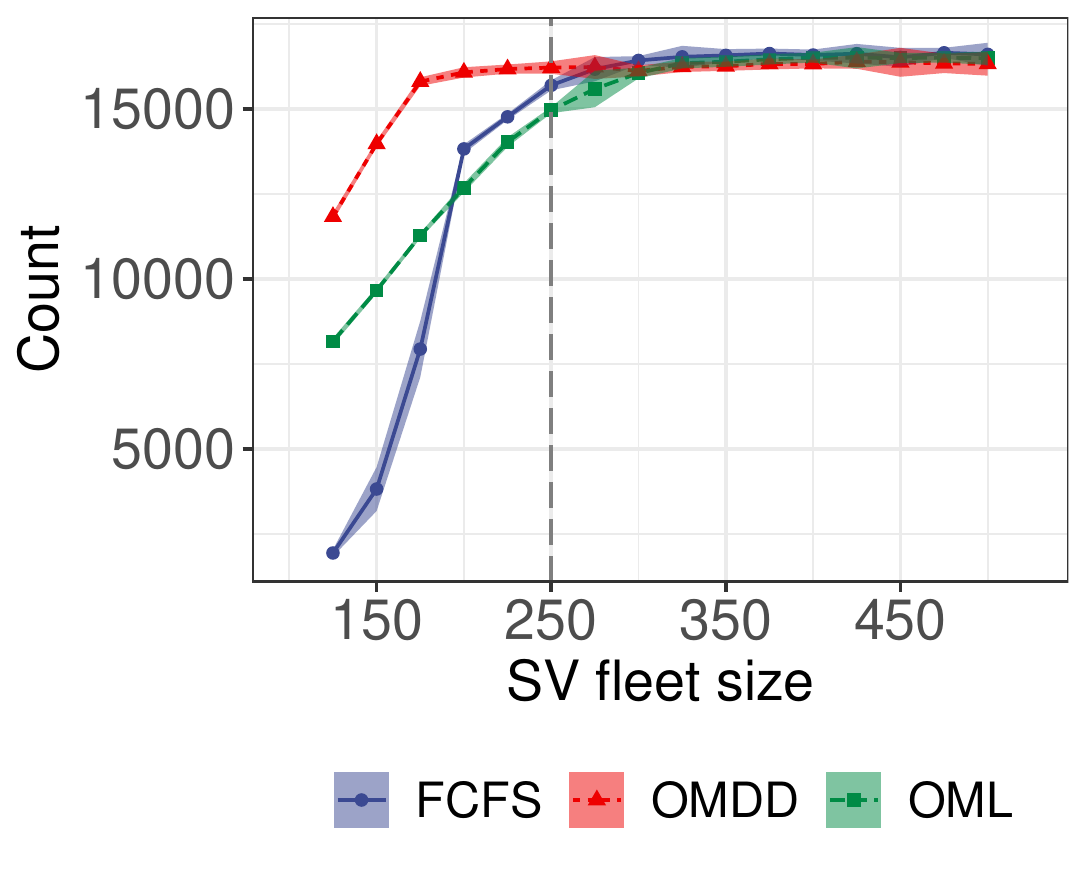}} \\
    \subfloat[Charging requests fulfillment rate~\label{fig:caas_dynamics_satisfy}]{\includegraphics[width=0.33\linewidth]{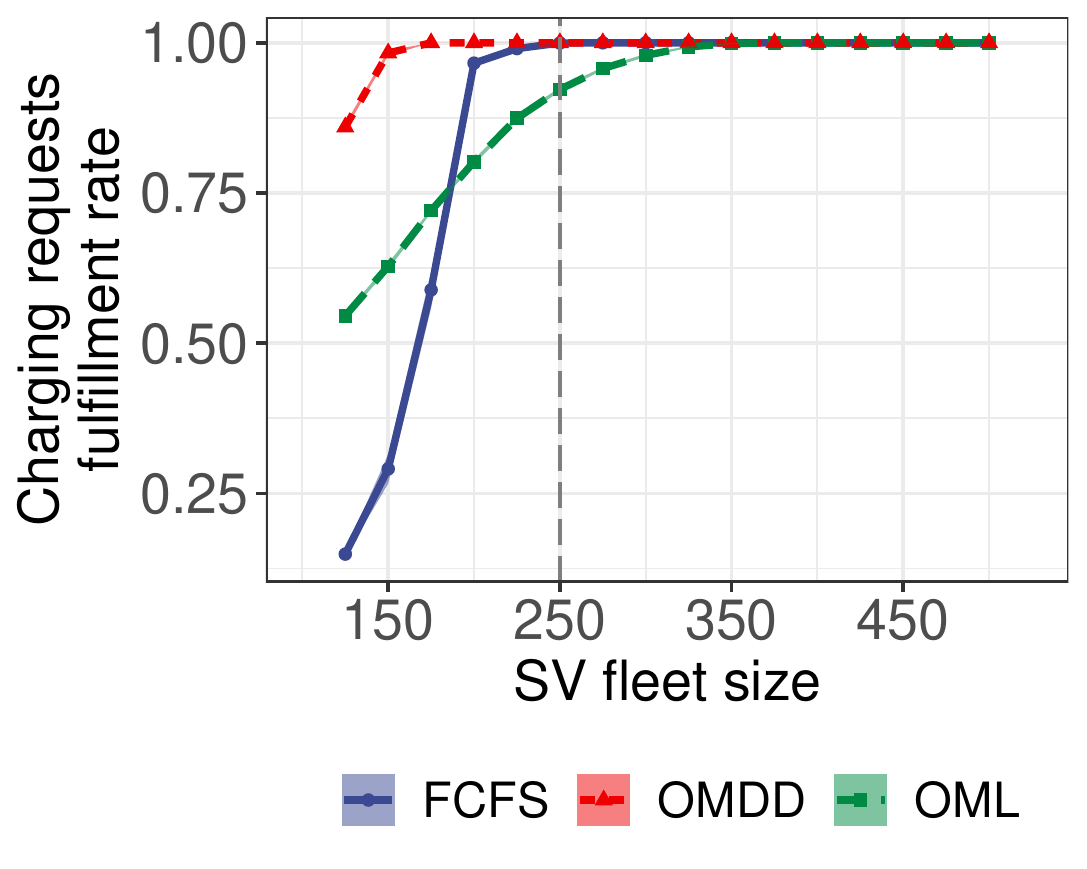}} 
    \subfloat[SV utilization rate~\label{fig:caas_dynamics_untilized}]{\includegraphics[width=0.33\linewidth]{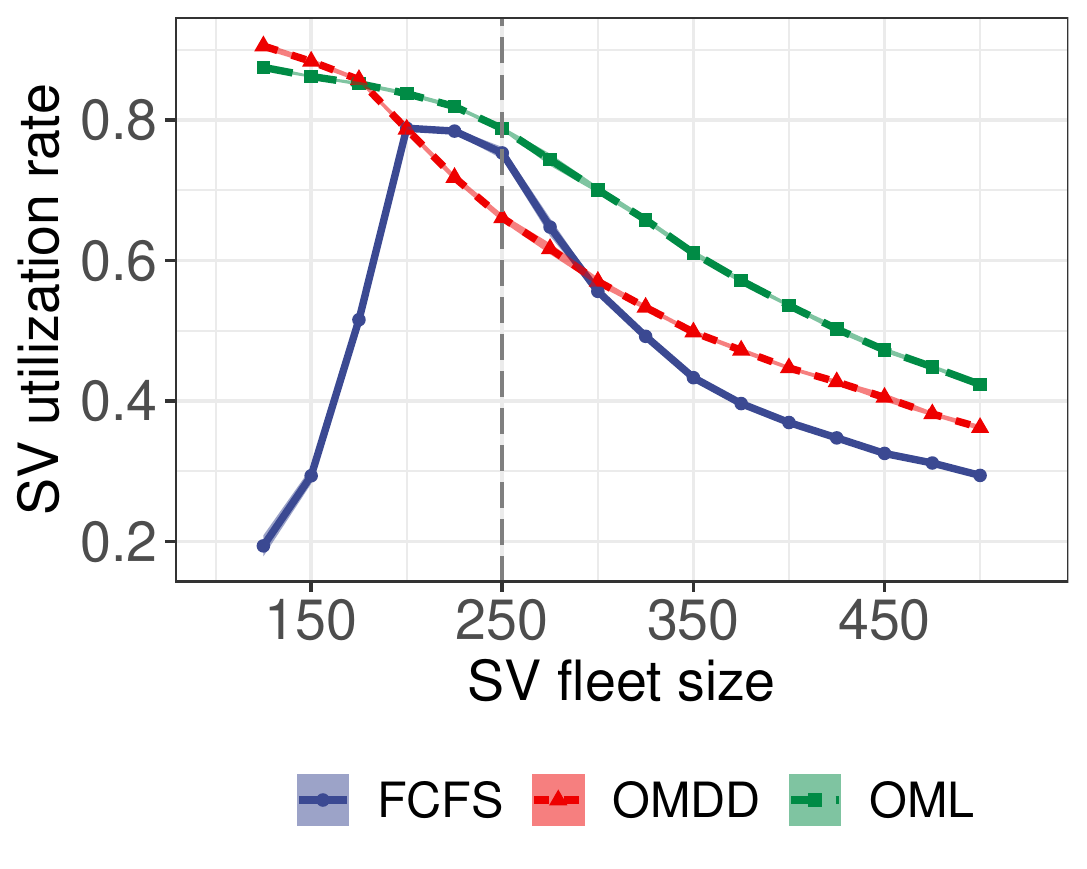}} 
    \subfloat[Trade-off between fulfillment and utilization~\label{fig:caas_dynamics_U_S}]{\includegraphics[width=0.33\linewidth]{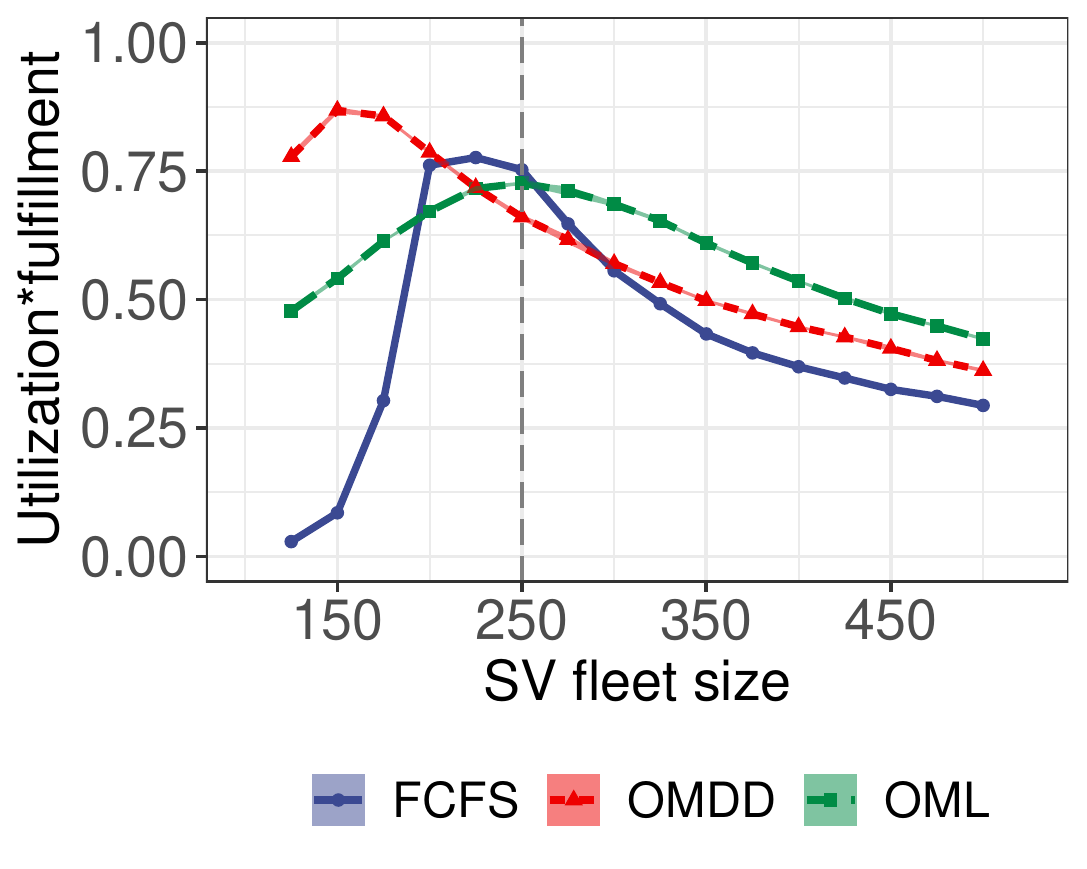}} \\
    \caption{CaaS dynamics under the SV fleet size of 125 to 500, incremented by 25}
    \label{fig:caas_dynamics}
\end{figure}

We first calibrate the impacts of various SV fleet size based on four major performance metrics as shown in Figure~\ref{fig:caas_dynamics}, including the average service waiting time, the SV fleet utilization rate, the fulfillment rate of charging requests, and the trade-off between service fulfillment rate and the fleet utilization. The results shown are the average system performances from 5 AM to 12 AM, and the shaded areas present the variation of system performances over simulation runs. The performances of all three dispatching strategies are evaluated under the same set of random seeds. It can be directly observed that the performances of the CaaS are stable at the system level with minimum deviation under stochastic demand. Figure~\ref{fig:caas_dynamics_wait} presents the average waiting time for all served LDEVs on the logarithmic scale, with SV fleet size varying from 125 to 500. The FCFS policy is found to be inferior to the other two policies due to its myopic dispatching strategy. On the other hand, the decreases in average waiting time between the OML and the OMDD policies are found to differ with the increasing number of SVs. The waiting time of both policies plateaus at around 1.5 min with a sufficient number of SVs (400 and over), indicating the potential of prompt responses from both dispatching strategies to fulfill the charging requests efficiently. Unlike the OML policy, the average waiting time for OMDD is decreased slowly between $N_{SV}$=225 and $N_{SV}$=325 while dropping faster at the other two ends. The reason is due to the fundamental differences in dispatching SV between the OML and OMDD policies. When there is an undersupply of SVs, the OMDD tends to serve more charging requests than the OML policy as it conducts the matching for the entire SV fleet rather than the available SVs. In many cases, this will sacrifice the waiting time of nearby requests to meet the requests that are distant from the SV fleet. As can be seen in Figures~\ref{fig:caas_dynamics_satisfy_count} and~\ref{fig:caas_dynamics_satisfy}, the OMDD policy is close to meet all charging needs with the fleet size of 225, where fulfilling further requests will lead to more frequent long-distance relocations hence resulting in a slower reduction in average waiting time. This issue is later resolved when there is an oversupply of SVs (325 and over) for the OMDD policy.

In addition to the performances on the demand side, we also evaluate the performances from the supply side by calculating the SV utilization rate. The rate measures the efficiency of the SV fleet based on the amount of time they spent in R and S states. With a small fleet size ($<$200), the SVs are found to spend over 80\% of the time serving charging requests, with 10\% of the time in visiting the depot for reloading MBUs. On the other hand, the FCFS policy can only make use of less than 70\% of the fleet resources since many requests are without reach due to the inefficient SV dispatching under the greedy scheme. While more SVs are found to reduce the average waiting time, there is also a notable drop in the fleet utilization rate. This implies the waste of resources, where on average, over 50\% of the SV service time is in the idle state if the fleet size is greater than 350 for the OMDD policy and 450 for the OML policy. And the differences in utilization rate under the same fleet size suggests that the OMDD policy is more efficient than the OML policy in satisfying the charging requests while delivering a comparable level of service for the LDEVs. 

Based on the previous discussion, it is evident that the CaaS system suffers the same unbalanced supply and demand issue as in the MaaS industry. And it is crucial to identify a cost-effective fleet size that seeks a compromise between the level of service and the utilization of supply resources. In this regard, we can examine how the product of charging requests fulfillment rate and the fleet utilization rate may change with increasing fleet size so that we can identify the desirable fleet size for the CaaS system. The results are summarized in Figure~\ref{fig:caas_dynamics_U_S}, with the product of 1 representing the ideal system performance and 0 for the worst case. The figure suggests an optimal fleet size of around 150 for the OMDD, 225 for the FCFS, and 250 for the OML policy. This observation again supports the superiority of the OMDD policy over the other two alternatives in empowering an efficient CaaS system. Finally, we choose the SV fleet size of 250 to explore further details on the operation dynamics of the CaaS system, which is found to yield a high request fulfillment rate, strike a balance between off-peak and peak hour charging demand and achieve good overall system efficiency for all three dispatching policies. 

\subsubsection{Spatial and temporal dynamics}
\begin{figure}[H]
\centering
    \subfloat[Hourly average waiting time~\label{fig:avg_wait_time}]{\includegraphics[height=.3\linewidth,width=0.4\linewidth]{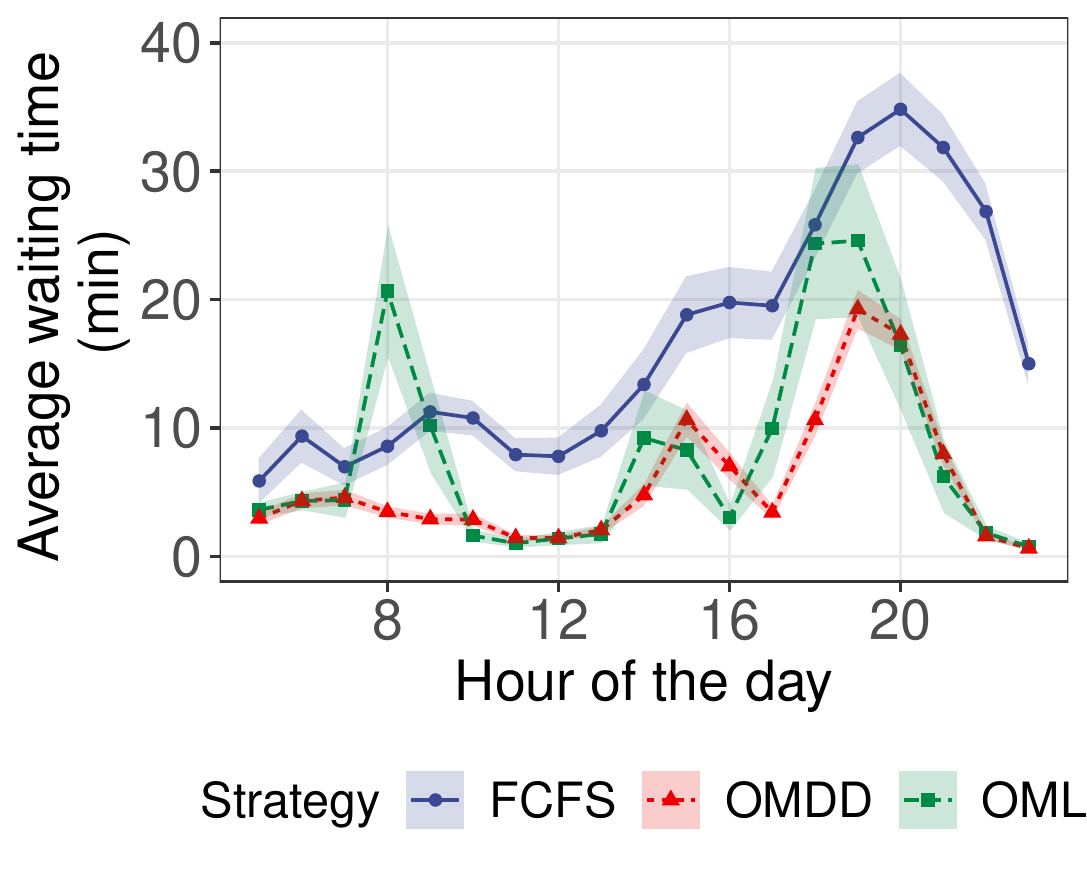}}
     \subfloat[Hourly average relocation distance~\label{fig:avg_relo_duration}]{\includegraphics[height=.3\linewidth,width=0.4\linewidth]{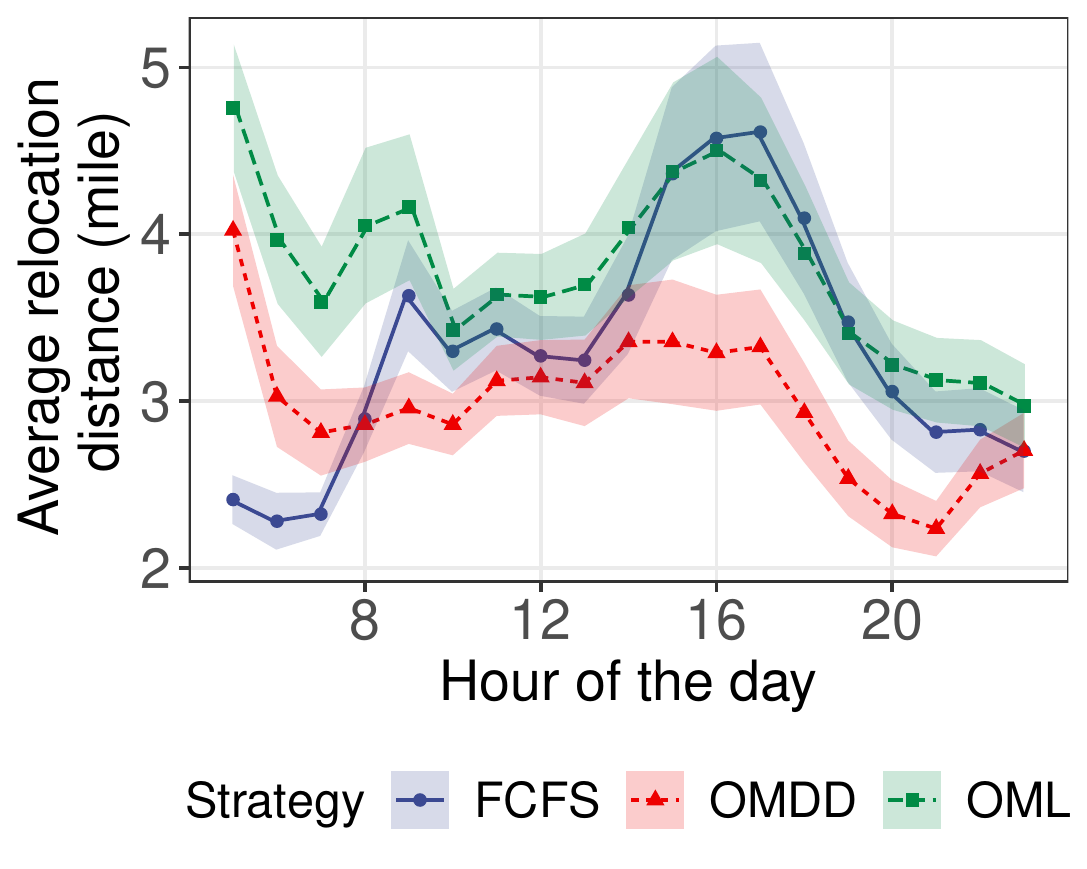}}\\
    \subfloat[Total number of generated and served LDEV charging requests (15 min interval)~\label{fig:yy_served_EV_demand}]{\includegraphics[width=.8\linewidth]{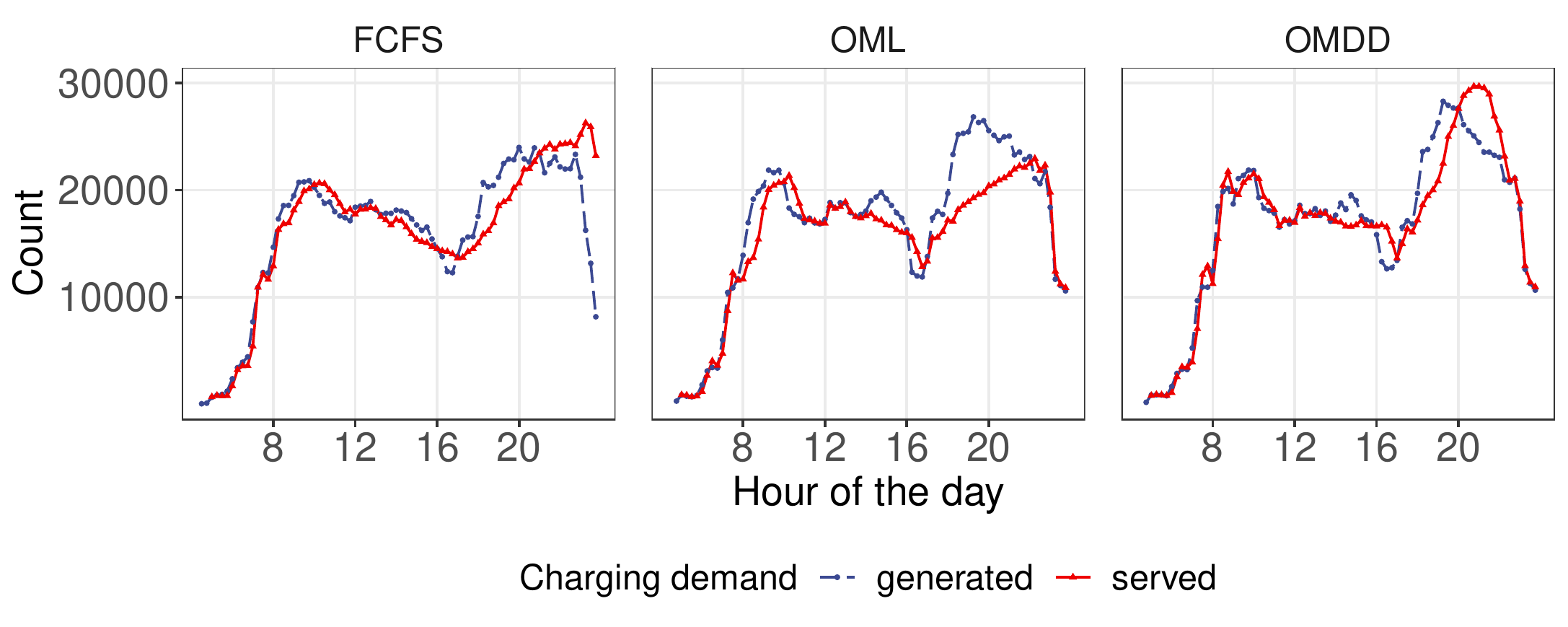}}
\caption{Temporal variation of average waiting time and relocation distance}
    \label{fig:temporal_analyses_served}
\end{figure}

\begin{figure}[H]
\centering
    \subfloat[Average waiting time~\label{fig:map_avg_wait}]{\includegraphics[width=.8\linewidth]{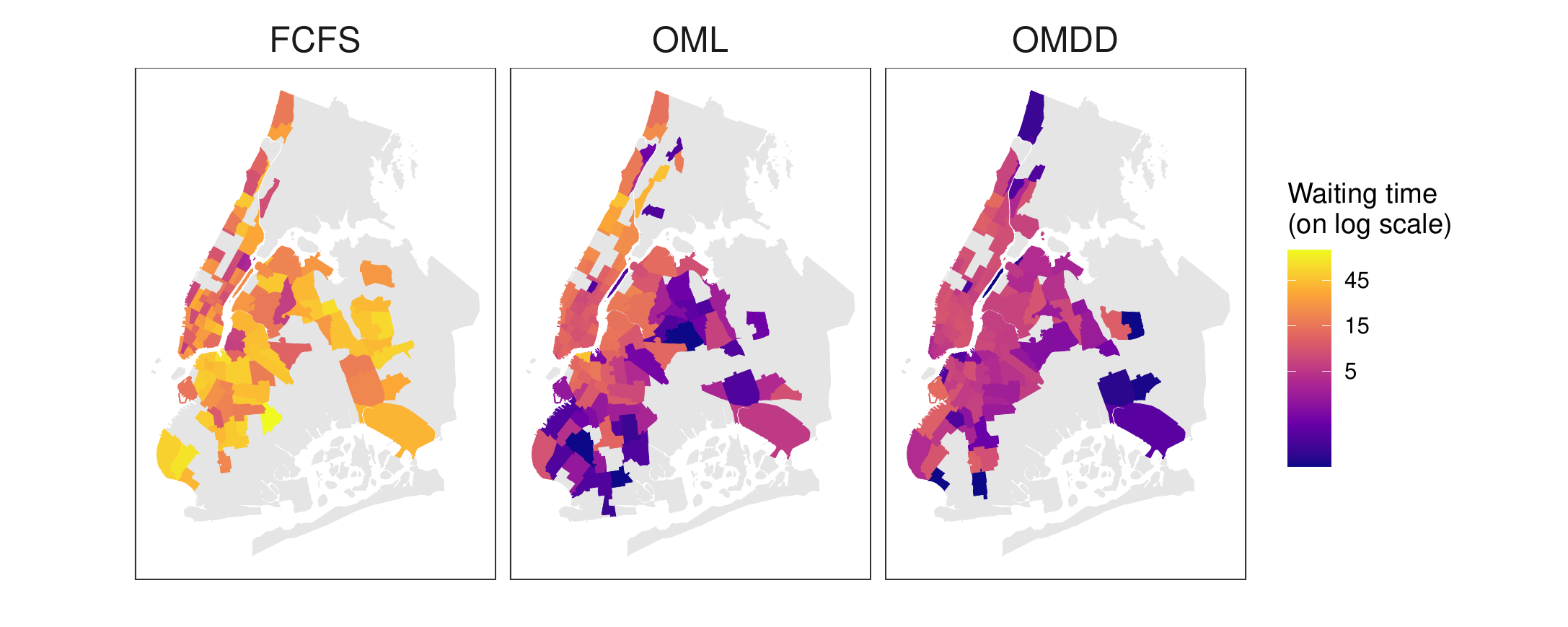}}\\
    \subfloat[OD distribution of SV relocation trips (Only OD pairs with more than 10 trips are displayed)~\label{fig:SV_od}]{\includegraphics[width=.8\linewidth]{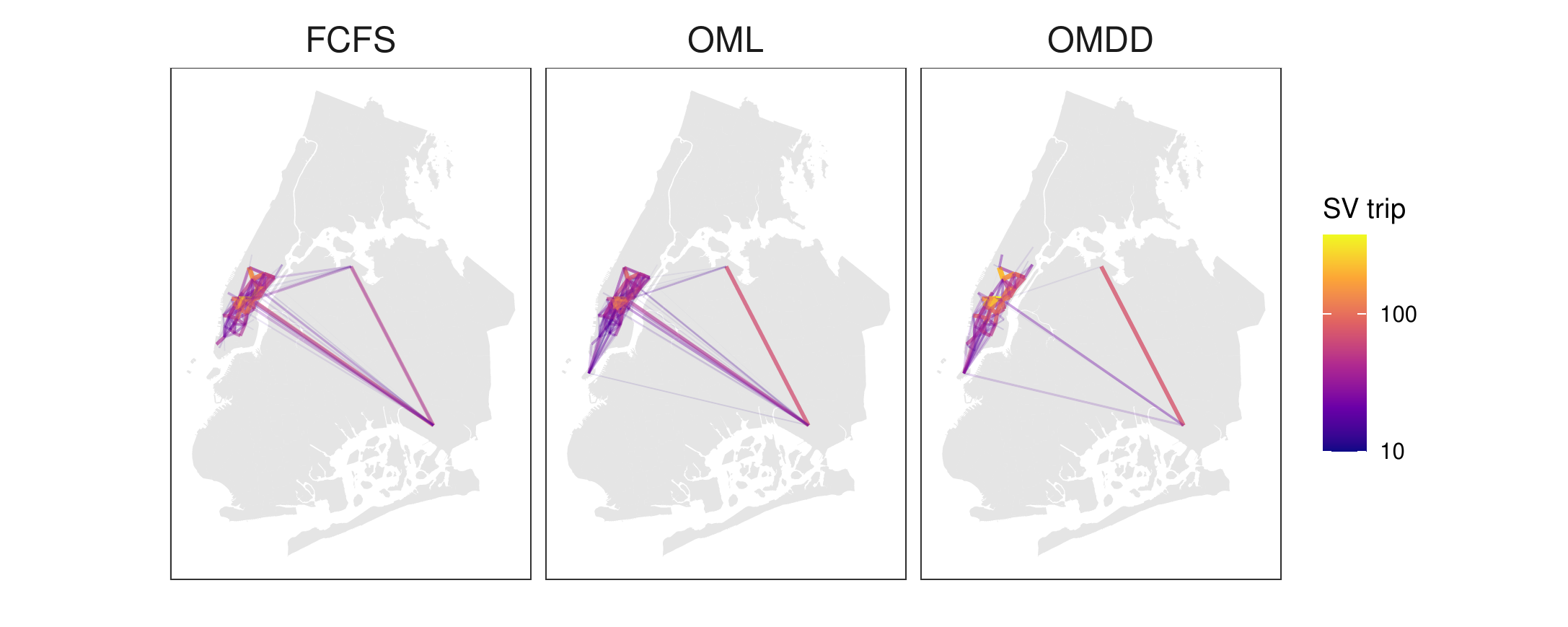}}\\
    \subfloat[SV fleet utilization rate percentage~\label{fig:map_work_N_all}]{\includegraphics[width=.8\linewidth]{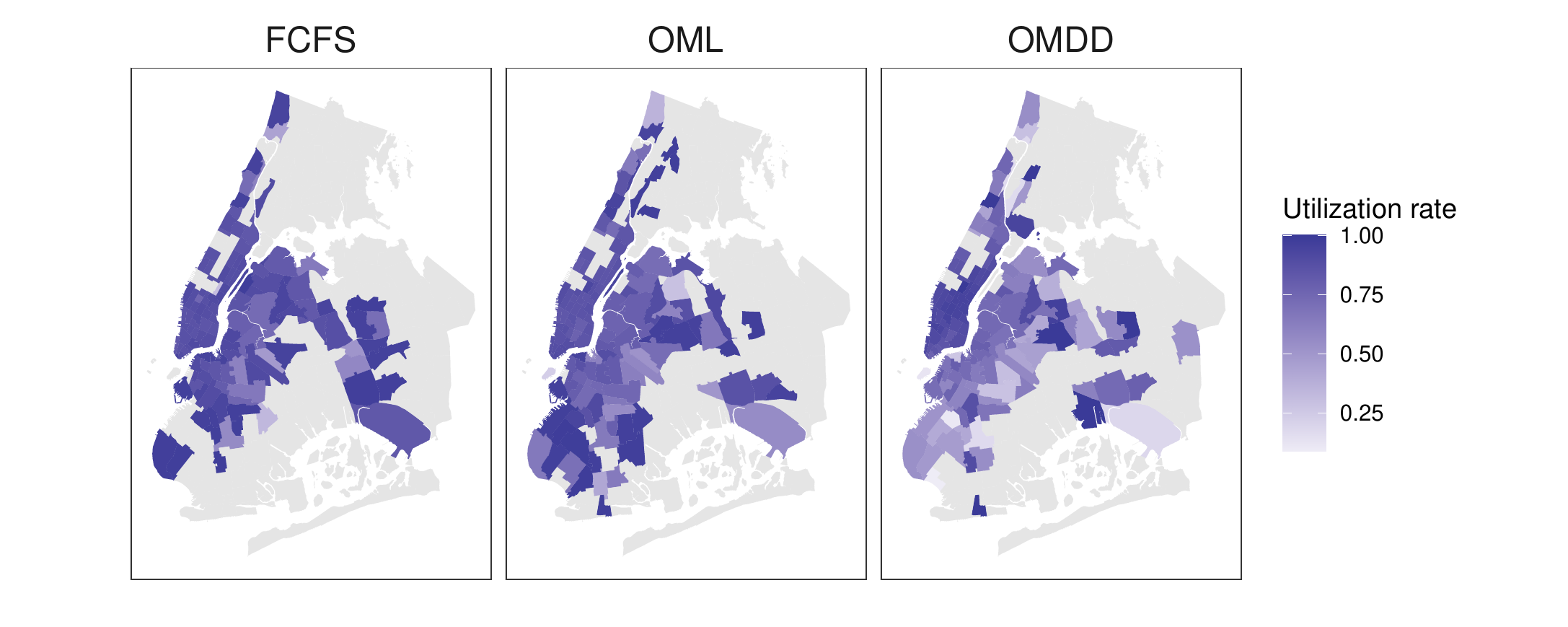}}
    \caption{Spatial distribution of the average waiting time}
    \label{fig:spatial_dynamics}
\end{figure}

Figure~\ref{fig:temporal_analyses_served} demonstrates the temporal dynamics of the average waiting time in the CaaS system from 5 AM to 12 AM, and Figure~\ref{fig:spatial_dynamics} shows the spatial distribution of system performances across the study area. For the hourly average waiting time in Figures~\ref{fig:avg_wait_time}, we observe three notable surges in the average waiting time at around 8 AM, 3 PM, and 7 PM, which are associated with the increase in charging requests as shown in Figure~\ref{fig:yy_served_EV_demand}. We note that the 3 PM demand peak for the OML and OMDD policies is not observed in the original charging demand distribution in Figure~\ref{fig:served_demand_temporal}. Instead, it is the result of the altered system dynamics with the delayed service of charging requests in previous time steps. For the FCFS policy, a steady increase in average waiting time emerges after the morning peak and continues until 8 PM, which is the result of the persistent rollover of unsatisfied charging requests caused by myopic dispatching decisions. The OML policy shows a significant disparity between peak and off-peak periods, with off-peak waiting time being less than 5 minutes and peak-hour waiting time exceeding 20 minutes. On the contrary, the OMDD policy is found to reach more consistent performances under the shorter surge in charging requests. And the performance during morning peak hours is comparable to the off-peak periods. The performances deteriorate during the evening peak due to higher demand level and longer duration but are still more resilient to demand changes than the OML policy with the same number of SVs. Similar observations can be identified based on the average waiting time's spatial distribution, as shown in Figure~\ref{fig:map_avg_wait}. The FCFS policy is found to be worse than the other two alternatives, with the average waiting time over 15 minutes in almost all taxi zones. The performances are better in high demand areas such as Manhattan but are considerably worse in other places. Both OMDD and OML policies are able to achieve comparable performances in and outside Manhattan because of the optimization performed at the system level. We also note that the average waiting time of the OMDD policy is notably better than the OML policy within high demand areas, and this can be explained by the fundamental differences of the dispatching philosophies as highlighted in Figure~\ref{fig:SV_od}. Specifically, the OMDD shows more clustered relocation trips within Manhattan and much fewer back and forth relocation trips between Manhattan and other areas such as the JFK airport than the other two policies. This attributes to the implemented dynamic detour strategy, which reduces relocation distance by inserting nearby requests between existing charging and often divides one long relocation trip into multiple short segments to fulfill more charging requests. Meanwhile, the OML policy makes one-shot relocation decisions for the SVs that frequently lead to relocation trips between high and low demand areas, with an example of more OD connections between Manhattan and the JFK airport in Figure~\ref{fig:SV_od}. As a consequence, the SV fleet utilization rates of the FCFS and the OML policies are much higher than the OMDD in taxi zones outside Manhattan, as can be seen in Figure~\ref{fig:map_work_N_all}). In the majority of the cases, the SVs in these zones will be immediately sent back to Manhattan under the OML and FCFS policies under the one-step decision making. But the OMDD policy tends to set aside these vehicles and look for alternate SVs in other zones through the multi-step optimal relocation strategy. This helps to avoid a significant number of unnecessary relocation trips, with a notable example of the much lower utilization rate at JFK airport under the OMDD policy. 

The advantage of the OMDD policy is also captured by the changes in the average relocation distance, as shown in Figure~\ref{fig:avg_relo_duration}. Before 8 AM, the average relocation distance starts with high initial values for the OMDD and OML policies (4 and 4.7 miles) as the charging requests are more sparsely distributed spatially, hence requiring longer relocation trips and exhibiting fewer opportunities for the dynamic detour of the OMDD policy. The FCFS, on the other hand, starts with a short average distance as a result of assigning the nearest available SVs to the charging requests under relatively low charging demand. And its relocation distance starts to increase and remain at a higher level during the day due to insufficient supply, which may no longer sustain the short distance even for the nearest available SVs. The OML policy also shows a similar relocation distance as the FCFS policy. The reason is because of more frequent back and forth relocation for immediate available SVs. However, the OMDD policy provides a better solution to mitigate the back and forth relocation issue with its dynamic detour feature and its ability to coordinate among the entire SV fleet. 

\begin{figure}[H]
    \centering
    \includegraphics[width=\linewidth]{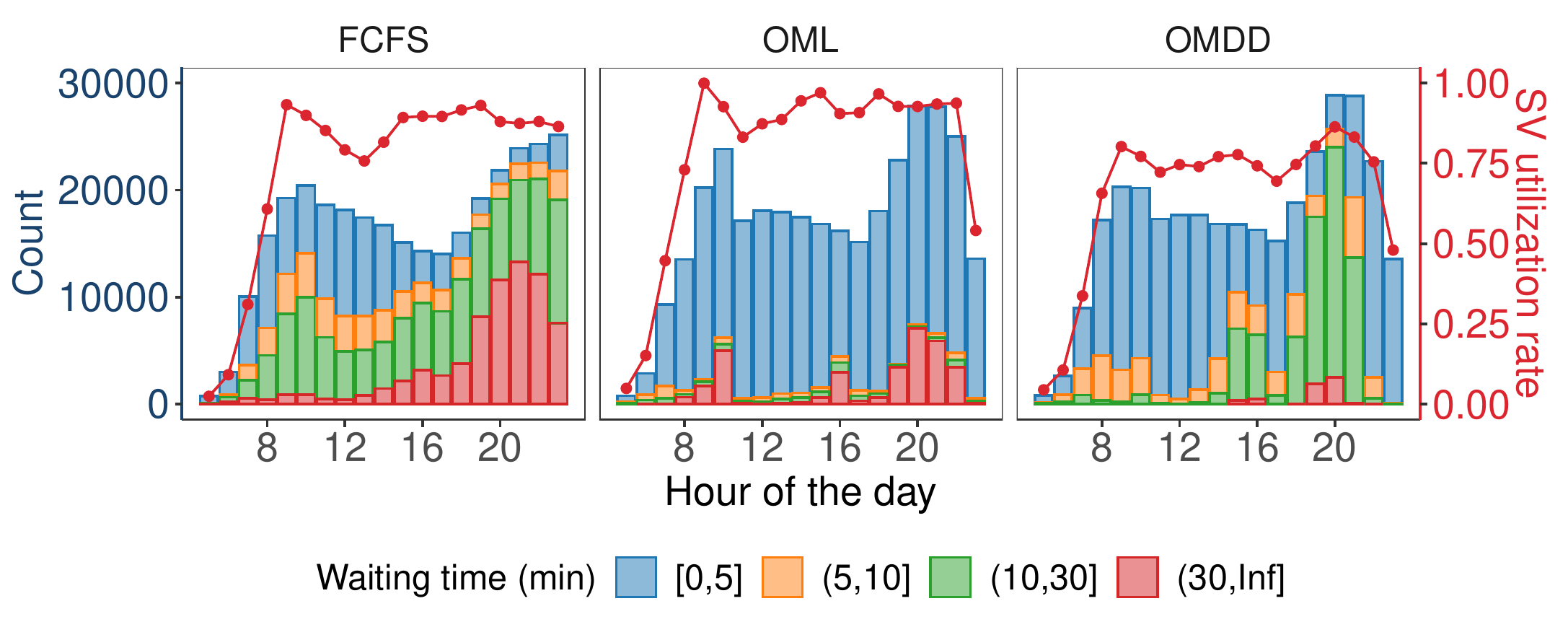} 
    \caption{Stacked histogram of charging demand by waiting time~\label{fig:wait_time_level}}
    \label{fig:supply_demand}
\end{figure}

To further demonstrate the performances of different policies during peak hours, we divide the waiting time of the charging requests into four groups and plot the stacked histogram along with the corresponding SV fleet utilization rate in Figure~\ref{fig:wait_time_level}. The distributions of waiting time are significantly different during peak hours for the OMDD and OML policies. The OML policy is found to serve a number of charging requests within 5 minutes but also keeps many charging requests waiting for more than 30 minutes. This is less of a problem during off-peak periods. But it raises a major concern when there is an undersupply of SVs in prime time, where excessive leftover requests are generated and lead to the snowballing of waiting time. This can also be confirmed by the distribution of served charging requests shown in Figure~\ref{fig:yy_served_EV_demand}. The OMDD, on the other hand, is found to adapt to demand surge at the morning peak (8- 10 AM) and minimize the number of rollover requests during the evening peak (7 - 10 PM). And this is the result of better utilization of the SV fleet during peak hours and the saving of supply resources before the arrival of the peak hours, as reflected by the change of SV utilization rate for the three policies.  

\subsection{Service fairness}
We are not only interested in the average metrics of the system performances but also the fairness of the CaaS regarding both demand and supply sides of the system. To this end, we introduce Gini coefficients~\cite{gini1912variabilita} to evaluate the spatial equity in terms of the differences of average waiting time over all taxi zones and the disparity of charging requests served by individual SVs. The calculation for the Gini coefficient is formulated as follows:

\begin{equation}
    G = \frac{\sum_i \sum_j \left| \bar{t}_i ^w - \bar{t}_j ^w\right|}{2 n^2 \bar{t^w}}
    \label{eq:gini}
\end{equation}

where $\bar{t}_i ^w$ is the average waiting time for zone $i$. $n$ is the total number of zones, and $\bar{t^w}$ is the mean value of the average waiting time of all zones. Similarly, we can perform equity analysis on served requests per SV by replacing $\bar{t}^w_i$ with the served charging requests for SV $i$.

\begin{figure}[H]
    \subfloat[Gini coefficients of average zonal waiting time~\label{fig:gini_wait_zone}]{\includegraphics[width=0.5\linewidth]{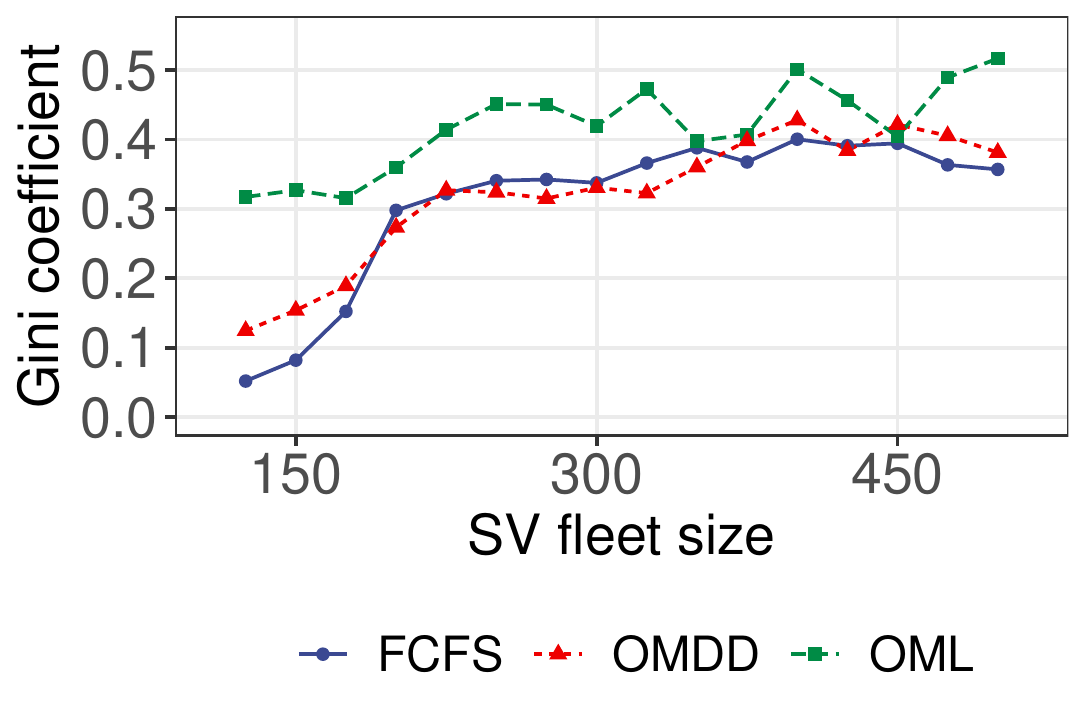}} 
    \subfloat[Gini coefficients of served charging requests per SV~\label{fig:gini_count_SV}]{\includegraphics[width=0.5\linewidth]{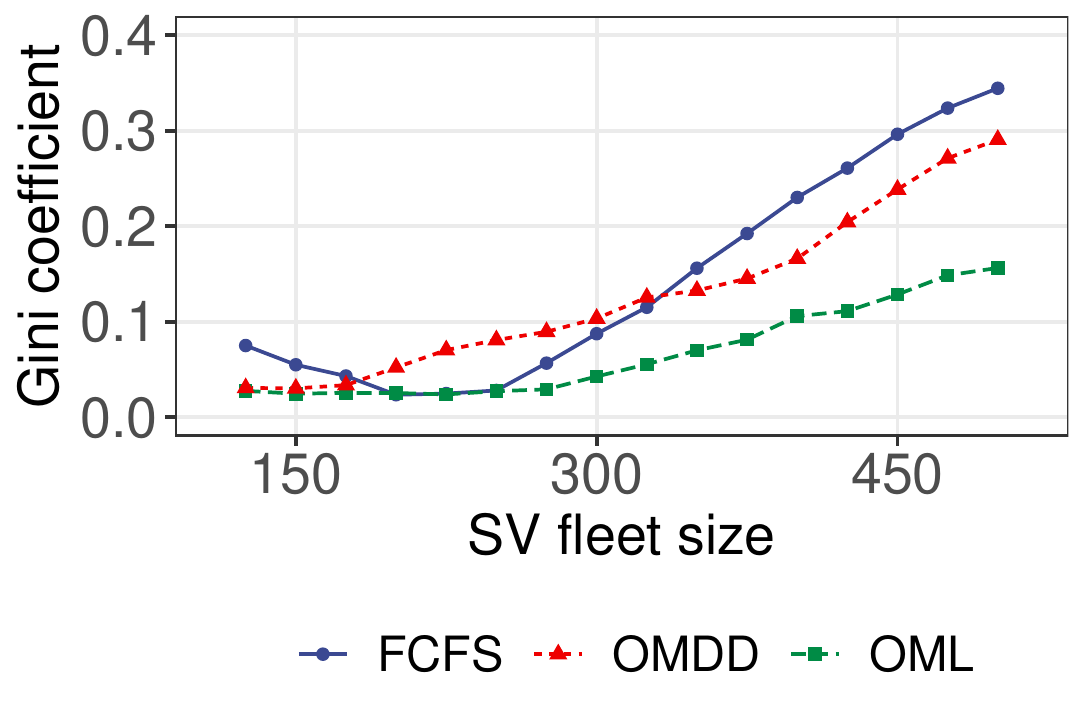}} \\
    \caption{CaaS performance on equity issues}
    \label{fig:gini_analyses}
\end{figure}

The results of the Gini coefficients are summarized in Figure~\ref{fig:gini_analyses}. For average waiting time, we report that the Gini coefficient increases with the SV fleet size (up to 500 SVs), and the Gini coefficient exceeds 0.3 with a fleet greater than 250 in all three policies. The increase of the Gini coefficient is due to the migration of the system state with all zones having long waiting time to the reduction of waiting time in certain zones, which is expected to decrease with further increase in SV fleet size as all zones receiving low waiting time. For the FCFS and OML policies, the surge of charging demand during peak hours motivates additional SVs to flow towards the high-demand areas, leading to an unbalanced SV fleet distribution. As a consequence, the waiting time in high-demand zones and zones close to the high-demand region is expected to be lower than in other areas on average. We note that the Gini coefficient of both FCFS and OMDD policies is lower than that of the OML policy. Nevertheless, the lower Gini coefficient for FCFS is because of the consistently worse performances over the entire study area, which is precisely the opposite of the reason for the lower Gini coefficient of the OMDD policy. 

As for the served requests per SV, the fairness of the assigned requests is vital for the CaaS system depending on the ownership of the SVs and how the SVs are operated. In the case of individually owned and operated SVs, the fairness of the dispatching policy is critical as it directly links to the SV drivers' workload and salary. As shown in Figure~\ref{fig:gini_count_SV}, we also observe that the Gini coefficient increases with more number of SVs. The initial low Gini coefficient is due to the lack of supply across all time periods, and all SVs are overloaded with charging requests. The FCFS policy is found to have the worst Gini coefficient with more SVs as a result of its nearest dispatching policy. Under the FCFS policy, SVs in the high demand tend to stay in these areas, which also applies to those in low demand areas, so that SVs in high demand areas will receive more requests with relatively shorter relocation trips. Despite the superior performances of the OMDD policy at the system level, it results in a consistently higher Gini coefficient than the OML policy, which is a side-effect of the dynamic detour strategy. Though the OML is making more relocation trips, it also ensures all SVs similar number of charging requests since it only focuses on immediately available SVs. On the other hand, as discussed in the previous section, the OMDD will generate clustered short trips in high demand areas, which will likely assign multiple requests with short relocation trips for some of the SVs while setting other SVs in the idle state. This represents a major barrier to implementing the OMDD policy when the SVs are owned by individual drivers, despite its superior system-level efficiency.  

\subsection{CaaS system benefits}
After discussing the operation dynamics, we next present the overall benefits of the CaaS system for the MaaS industry compared to the adoption of FCSs. And we focus on the CaaS system benefits in terms of the savings in the OoS time, the relocation distance, and the financial sustainability of the CaaS. 
\begin{figure}[H]
    \centering
    \subfloat[Savings in OoS time (min per charging request)\label{fig:system_oos_saving}]{\includegraphics[width=0.5\linewidth]{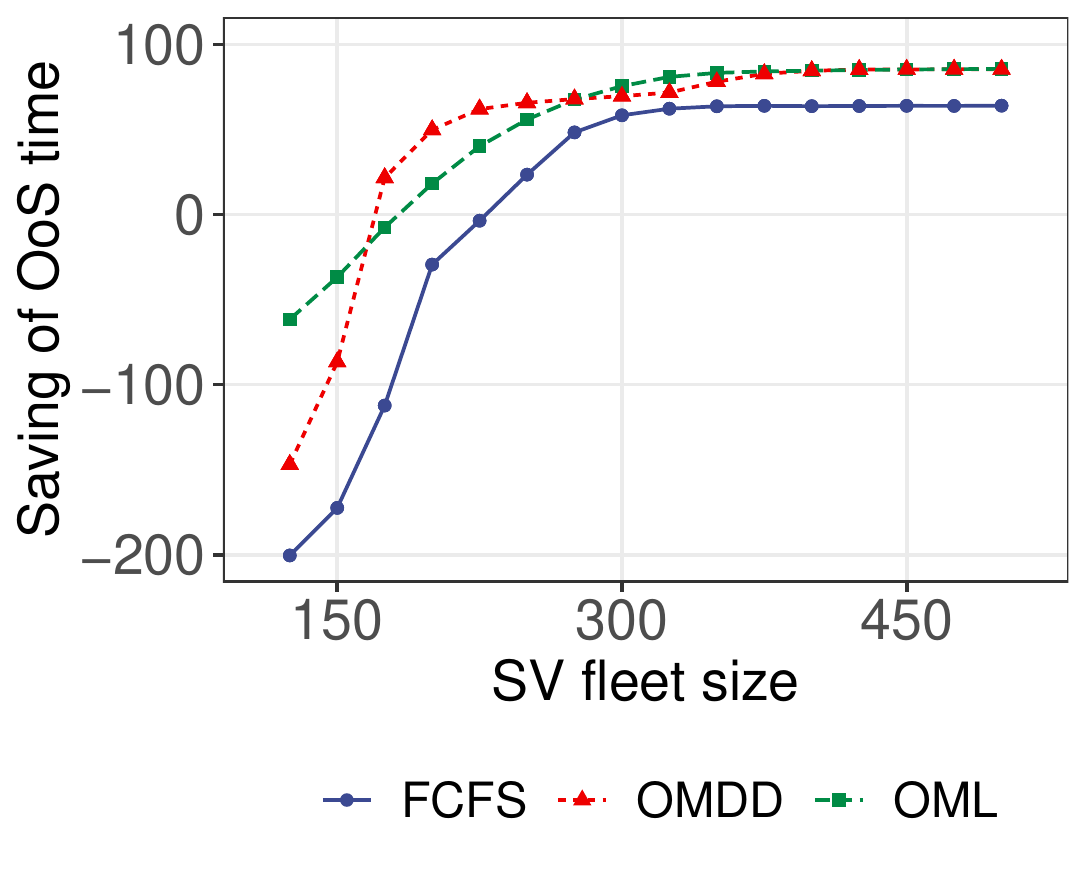}}\hfill
    \subfloat[System savings in relocation distance (mile per min)\label{fig:system_relo_saving}]{\includegraphics[width=0.5\linewidth]{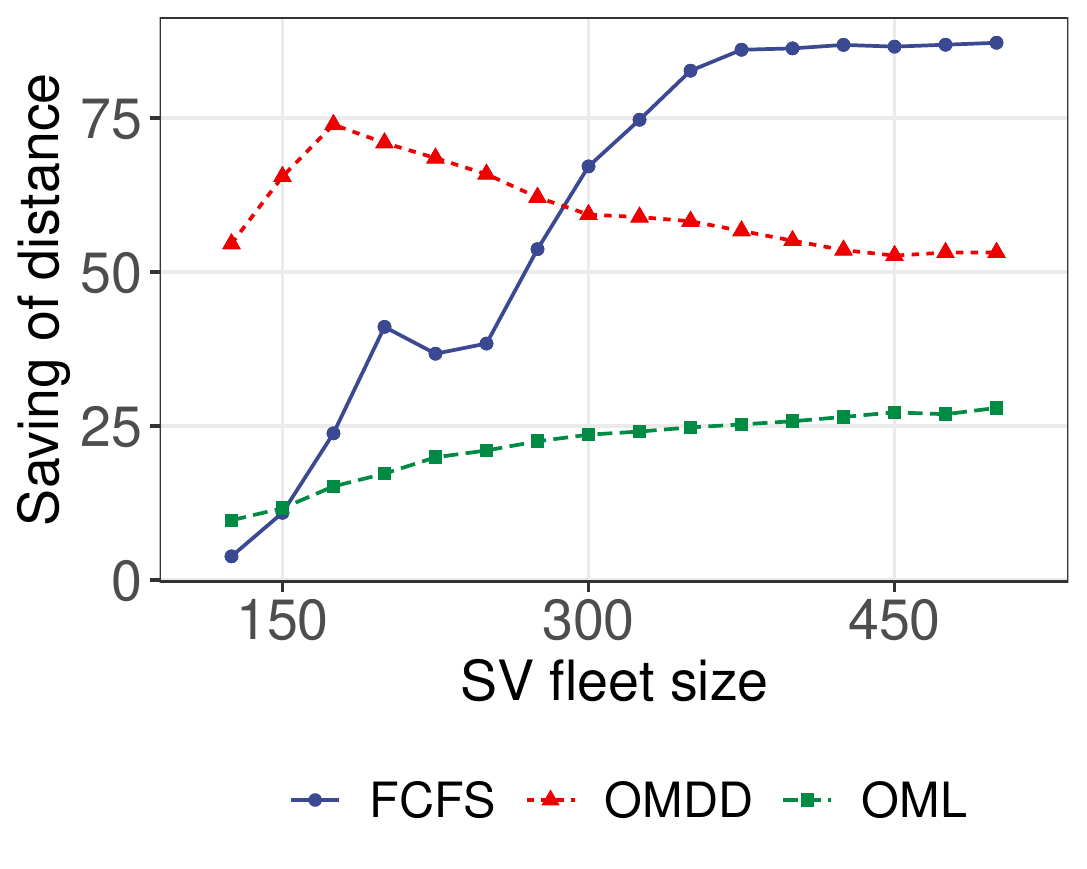}}\\
    \caption{CaaS system savings}
    \label{fig:savings}
\end{figure}
Figure~\ref{fig:system_oos_saving} indicates the savings of OoS time between the CaaS and charging through FCSs. All metrics from the CaaS are converted under the FCE consideration to make the results comparable under the same level of charging demand. The results show that the CaaS achieves consistent savings with a sufficient number of SVs ($>$250) for all three policies. In particular, the charging requests under the OML and OMDD policies receive similar and most considerable savings at over 4 minutes per charging request. Under the FCFS policy, the saving plateaus at 3 min with an SV fleet size larger than 300. The savings are lower than those under the OMDD and the OML, which is also evidenced in the hourly average waiting time in Figure~\ref{fig:caas_dynamics_wait}. The saving will lead to more served trips in the MaaS and is critical for promoting a wide adoption of EVs in the MaaS industry. 

Figure~\ref{fig:system_relo_saving} presents the system savings in travel distance due to the shift from charging at FCSs to the CaaS. The LDEV's relocation distance to the nearest FCS is assumed to be 15 miles on average. The saving of distance for charging under the OMDD policy remains at a high level of over 50 miles per minute, followed by the savings of the OML between 10 to 26 miles per minute. The savings of the OMDD policy can be translated into more than 5700 miles per day from 5 AM to 12 AM. The FCFS is found to achieve the highest saving in charging miles with a sufficient number of SVs, as the requests are assigned to the nearest SVs in sequence. For the OML, the saving also increases with more number of SVs since it increases the density of the SVs and naturally reduces the distance to reach the charging request. On the other hand, the OMDD policy reaches the highest saving at the fleet size of 175 but will have to trade the savings in charging distance for the reduction in average waiting time with a larger fleet. As a consequence, the CaaS also marks significant savings of energy in charging the large-scale LDEVs of the MaaS industry. 

\begin{figure}[H]
    \centering
    \includegraphics[width=\linewidth]{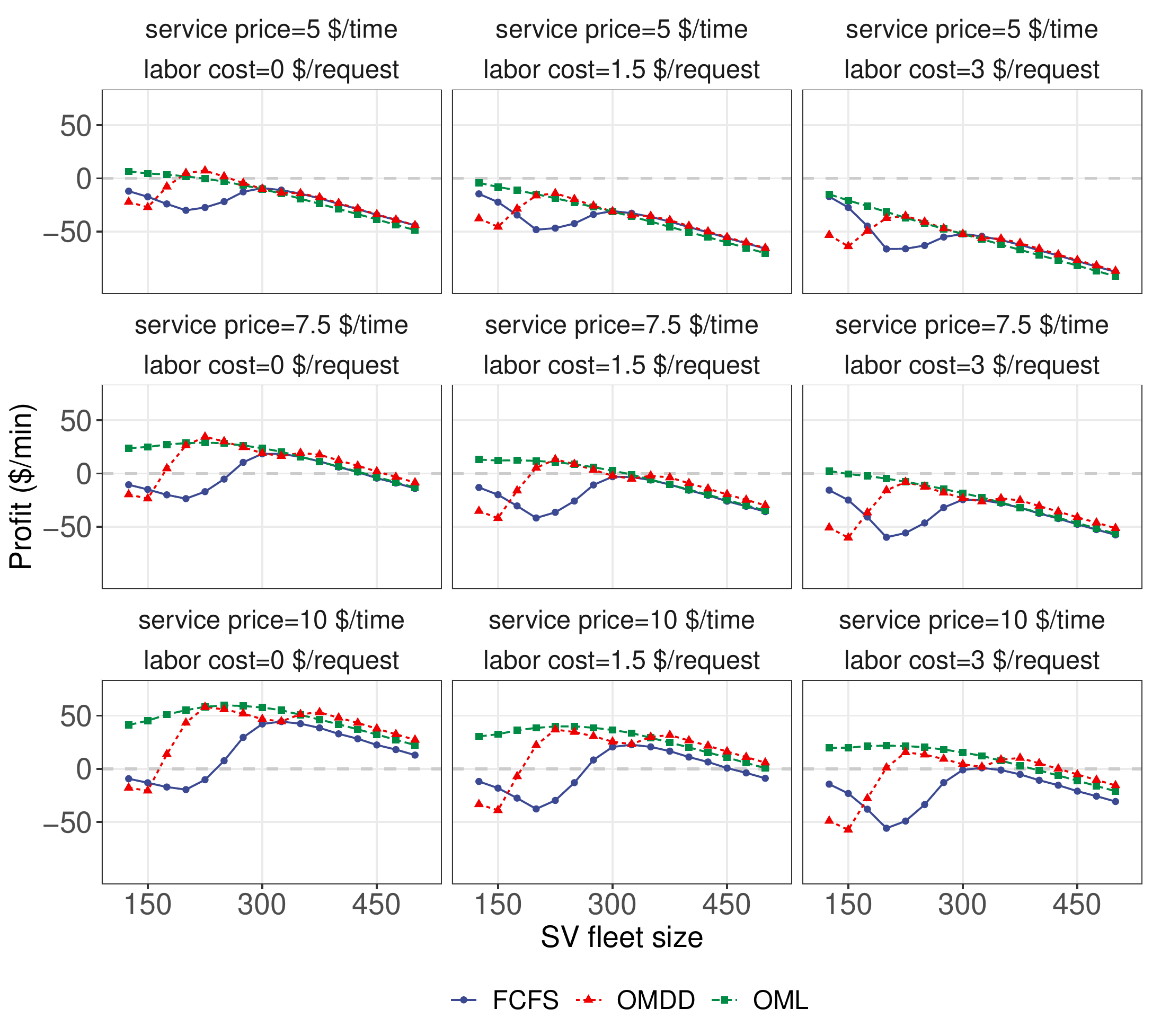}
    \caption{CaaS system profit ($\$$ per min)}
    \label{fig:profit}
\end{figure}

While demonstrating promising contributions to save OoS and charging distance, the feasibility of the CaaS depends highly on its financial sustainability. We conduct a sensitivity analysis on the profit level of CaaS with respect to the service price and the labor cost per charging request, and the results are shown in Figure~\ref{fig:profit}. In general, we find that the service price contributes significantly to the profit and CaaS users are willing to pay the extra to enjoy a higher level-of-service as compared to charging at FCSs. Besides, the labor cost also acts as a deterministic factor that affects the profit level, with the increase in labor cost resulting in a steep reduction of service profit. For the given service price and labor cost, the profit resembles a concave function with respect to the SV fleet size, and expanding the SV fleet for reduced average waiting time beyond an optimal fleet size will not square the additional capital and operation cost for the CaaS. And this optimal fleet size is found to shift to the left with increasing labor cost and shift to the right with higher service price per charging request. In this regard, for a CaaS startup with a limited initial budget, a higher salary paid to the SV drivers means fewer funds to enlarge the SV fleet size, which then translates into less saving of relocation time and higher waiting time of the served LDEVs. On the other hand, a reduced labor cost would benefit the CaaS to gain a higher budget for a larger SV fleet, achieve better service levels, and ensure sustainable cash flow with improved profit levels. A special case is the $\$0$ labor cost, which can be achieved through an automated CaaS system, and this case performs considerably better than other cases in terms of financial sustainability. And this result highlights the autonomous MBU delivery service as an important component for a sustainable CaaS system.

Additional insights into the CaaS operation dynamics are observed by inspecting the profit changes with different configurations of the SV fleet under three dispatching strategies. These trends of the profit can be interpreted by two main components, as indicated in equation~\ref{eq:profit}. The first component is the capital and operation cost corresponding to the increase in SV fleet size, and the other component is the revenue from the potential CaaS users based on the service level with more SVs introduced to the fleet. Despite performing worse than the OMDD policy at the system level, the OML policy outperforms the OMDD for a higher system profit with a smaller fleet size (e.g., $N_{SV}<350)$. In this case, the OML serves more LDEVs with a shorter waiting time, which leads to a higher turnout rate for using the CaaS and, therefore, more revenues. 

As the SV fleet size exceeds 400, the profit under the OMDD policy slightly outperforms the OMDD as it offers a similar service level at the cost of fewer relocation miles. Based on the profit analyses, we conclude that the CaaS is a potentially profitable business under appropriate service price, labor cost, fleet size, and dispatching policy configuration. Moreover, despite a loss of revenue, more SVs will yield a higher saving of OoS time, which may well justify the reduced profit with more served trips in the MaaS industry.

\section{Conclusion and Future work}

In this paper, we propose the HABM to simulate the CaaS system dynamics, where SV fleet are dispatched to provide on-demand MBU delivery for the LDEVs in the MaaS industry. And we investigate the applicability and effectiveness of the CaaS system under three different dispatching strategies. The performances of the CaaS framework are validated through comprehensive numerical experiments following the same trip dynamics in the NYC taxi market where all taxis are assumed to be replaced with the LDEVs. The CaaS system is found to achieve high level of service and satisfy over 95 \% charging requests from the LDEVs in NYC with an SV fleet size of 250. And the CaaS system has the potential to achieve significant savings in OoS time and charging miles and maintain its financial sustainability under appropriate fleet and price settings. 

Based on the analyses of the operation performances, system benefits, and financial sustainability, we summarize our recommendation of the CaaS system configurations considering different types of stakeholders as follows:

\begin{itemize}
    \item TNCs: TNCs such as Uber and Lyft are considered large-scale operators aiming to provide a high level of service to their LDEV fleet and prioritize short waiting time and savings in OoS. The OMDD policy is recommended for TNCs with their own SV fleet, while the OML policy is more applicable if the TNCs decide to operate the SVs through crowd-sourcing. In this regard, the applicable SV fleet size can be shifted to the right of the profit-maximization ones, which will gain slightly lower profits but achieves less waiting time and more savings of OoS duration. The loss of profit in the CaaS is likely to be compensated for through extra gains from their trip services via the LDEV fleet. 
    \item Small startups as independent CaaS providers: Small startups operate CaaS with a budget constraint will need to balance labor costs and the investment into SVs to get the most profits. In this case, the OML protocol can be adopted, and the optimal number of SVs may be directly applicable to small startups who run CaaS for maximizing revenue. 
    \item State agencies that focus on environmental impact: The priority for these agencies is to reduce emission and energy consumption. They tend to discourage additional travel distances from the CaaS sector. In this study, the savings of travel distance under the OMDD doubles those under the OML. As a result, the agencies may consider subsidizing small startups to run the OMDD strategy or subsidize potential CaaS users for accepting a higher average waiting time.
\end{itemize}

While the study represents an initial attempt to investigate the dynamics of the CaaS system, it opens up several future research opportunities to advance our understanding of an optimally deployed CaaS system. First, the multi-depot setting can be investigated to improve the service level. Also, incentives and dynamic service price may motivate additional SVs to get out of the high-density area (e.g., Manhattan) to serve the LDEVs in other areas during peak hours so as to improve the unbalanced level of service in the city. In addition to the MaaS industry, the proposed CaaS system can also be extended to other micro-mobility sectors as a battery delivery service, e.g., the shared electric scooters, the electric bikes, and the electric commercial delivery drones. These scenarios describe a lite MBU delivery service under more uncertain trip trajectories and will require additional considerations to tailor to the unique features of these systems. 

\bibliographystyle{unsrt_new}

\bibliography{references}

\end{document}